\documentclass[a4paper,10pt]{article}
\usepackage{amssymb}
\usepackage{color}
\usepackage{graphicx}
\usepackage{amsmath}
\usepackage{amsmath, upgreek}
\usepackage{amsmath,tabu}

\newcommand{\C}{\mathbb{C}}
\newcommand{\R}{\mathbb{R}}
\newcommand{\Z}{\mathbb{Z}}
\newcommand{\N}{\mathbb{N}}
\newcommand{\cL}{\mathcal{L}}

\newcommand{\cuad}{{\sqcap\kern-.68em\sqcup}}

\newcommand{\norm}[1]{\|#1\|}

\newcommand{\s}{s}
 
 \newcommand{\cI}{\mathcal{I}}
  
 \newcommand{\cA}{\mathcal{A}}

\newcommand{\cB}{\mathcal{B}}

\newcommand{\cE}{\mathcal{E}}
\newcommand{\cH}{\mathcal{H}}
\newcommand{\cK}{\mathcal{K}}
\newcommand{\cO}{\mathcal{O}}
\newcommand{\cF}{\mathcal{F}}

\newcommand{\cP}{\mathcal{P}}
\newcommand{\vs}{\smallskip}
\newcommand{\cG}{\mathcal{G}}
\newcommand{\cM}{\mathcal{M}}
\newcommand{\cQ}{\mathcal{Q}}

\newtheorem{theorem}{Theorem}[section]
\newtheorem{proposition}{Proposition}[section]

\newtheorem{lemma}{Lemma}[section]
\newtheorem{corollary}{Corollary}[section]
\newtheorem{remark}{Remark}[section]
\newcommand{\bremark}{\begin{remark} \em}
\newcommand{\eremark}{\end{remark} }

 \renewcommand{\k}{\mathbf{k}}
\renewcommand{\j}{{\mathbf{j}}}
\DeclareMathOperator{\dist}{\rm dist}
\newcommand{\eps}{{\varepsilon}}
\newcommand{\ms}{\medskip}

\usepackage{cite}
\begin{document}

\begin{center}{\bf  \large   On $m$-order logarithmic Laplacians and related  propeties

 }

 \bigskip\bigskip

  {\small  Huyuan Chen\footnote{chenhuyuan@yeah.net}  }
 \bigskip

 {\small   Department of Mathematics, Jiangxi Normal University,\\
 Nanchang, Jiangxi 330022, PR China }  \\[4mm]

 The University of Sydney, School of Mathematics and Statistics, NSW 2006, Australia\\[18pt]

 \begin{abstract}

In this article,  we study $m$-order logarithmic Laplacian $\cL_m$, which is a singular  integro-differential  operator   with symbol $\big(2\ln |\cdot|\big)^m$ by the Fourier transform.
With help of these logarithmic Laplacians,  we build  the  $n$-th order Taylor expansion for fractional Laplacian with respect to the order
and the Riesz operators:  for $u \in C^\infty_c(\R^N)$ and  $x \in \R^N$,
$$
 (-\Delta)^s  u (x) =  u(x) + \sum^{n}_{m=1} \frac{s^m}{m!}\cL_mu(x) + o(\s^n) \quad  {\rm as}\ \,  \s\to 0^+ 
$$
and
$$
\big(\Phi_\s\ast u\big)(x) =  u(x) + \sum^{n}_{m=1}(-1)^m\frac{s^m}{m!}\cL_mu(x) + o(\s^n) \quad  {\rm as}\ \,  \s\to 0^+, 
$$
where $ (-\Delta)^s$ is the $s$-fractional Laplacian, $\Phi_\s\ast u$ is $s$-order of Riesz operator  with the form $\Phi_\s(x)=\kappa_{N,s}|x|^{2s-N}$ in $\R^N\setminus\{0\}$.  

  Moreover,   we analyze qualitative properties of these operators based on the order $m$,  such as basic regularity  and the Dirichlet eigenvalues.

 \end{abstract}

  \end{center}
  
  \noindent {\small {\bf Keywords}:    The $m$-order logarithmic Laplacian; Riesz kernel; fractional Laplacian.}\vspace{1mm}

\noindent {\small {\bf MSC2010}:     35R09,  35C20. }

\vspace{1mm}


\setcounter{equation}{0}
\section{Introduction}

 \qquad It is well-known  in \cite[p194, (2)]{GS} (see the Fourier transform of $r^{ \lambda}$) that the Taylor's and Laurent series expansion of $r^{ \lambda}$ with the formula  
\begin{equation}\label{1.1-series}
r^{ \lambda} =1+   \sum_{n=1}^m\frac{\lambda^n}{n!} ( \ln r)^n+o(\lambda^m) \quad  {\rm as}\ \,  \lambda \to 0^+.
\end{equation}
  Our aim in this article is to find a sequence of the integro-differential operators whose Fourier transforms having the form $(\ln r)^n$ with $n\in \N$, which gives the related series expansion  for the Riesz convolutions and the fractional Laplacians, their Fourier symbols $|\zeta|^{\pm2s}$ respectively,  as the order tends to zero, i.e.  $s \to 0$.  \vs

Let the integer $N\geq1$, $\s\in(0, \frac N2)$ and  
\begin{equation}\label{1.1-funda}
\Phi_\s(x)=\kappa_{N,s} |x|^{2s-N} \quad {\rm for}\ x\in  \R^N\setminus\{0\},
\end{equation}
where  $\kappa_{N,s}$ is a normalized constant with the form 
 \begin{equation}\label{1.1-kappa}
\kappa_{N,s} = 2^{-2s}\pi^{-\frac N2}  \frac{\Gamma(\frac{N-2s}{2})}{\Gamma(1+s)} s.
\end{equation}
Here $\Gamma$  is the  Gamma function,  both $\mathcal{F}$ and $\widehat \cdot$ denote the Fourier transform. 
With this constant, there holds in the distributional sense
$$\cF(\Phi_\s\ast u)(\zeta)=|\zeta|^{-2s}\hat{u}\quad {\rm for\ suitable}\ \, u.$$ 
The study of the Riesz kernel is an ancient area, known as the potential theory \cite{S,H}, double phase problems \cite{BY}, overdetermined problem \cite{LZ},
isoperimetric inequalities \cite{RRS} and more topics could see
\cite{CPV,D,HXY}.

 For $s\in(0,1)$, the function $\Phi_\s$ could also be viewed as the fundamental solution of the $\s$-order fractional Laplacian $(-\Delta)^\s$ with $\s\in(0,1)$   defined as a singular integral operator 
 \begin{equation}\label{fl 1}
 (-\Delta)^\s  u(x)=c_{N,\s} \lim_{\epsilon\to0^+} \int_{\R^N\setminus B_\epsilon(x) }\frac{u(x)-
u(y)}{|x-y|^{N+2\s}}  dy,
\end{equation}
where   $c_{N,\s}$ is a normalized constant 
$$c_{N,\s}=2^{2\s}\pi^{-\frac N2}\s\frac{\Gamma(\frac{N+2\s}2)}{\Gamma(1-\s)}.$$
With this normalized constant, there holds that for $u \in C^\infty_c(\R^N)$
\begin{equation}
  \label{eq:Fourier-representation}
\mathcal{F}((-\Delta)^\s u)(\xi) = |\xi|^{2\s}\hat u (\xi)\quad  {\rm for\ all} \ \xi \in \R^N.
\end{equation}

In recent years, integro-differential operators have been studied in various aspects  both by important applications and seminal advances in the understanding of nonlocal phenomena from a PDE point of view,  see e.g.
\cite{BCH,CS0,CS1,DT,KM,MN,RS,S07,S} and the references therein.  The fractional Laplacian is one of the most typical operator among nonlocal operators and it attracts the most attention.   It is well known that the fractional Laplacian has the following limiting
properties when $\s$ approaches the values zero and $1$:
 \begin{equation} \label{eq:limit-behaviour}
\lim_{\s\to1^-}(-\Delta)^\s u(x)=-\Delta u(x)
\quad\  {\rm and}\quad\ \lim_{\s\to0^+}  (-\Delta)^\s  u(x) = u(x)\quad\  {\rm for} \ \, u\in C^2_c(\R^N),
 \end{equation}
see e.g. \cite{EGE}.     A challenging problem is whether the expansion (\ref{eq:limit-behaviour}) could be extended to the   higher orders.  Also another  challenging topic is  what the expansion of the convolution $\s\mapsto\Phi_\s\ast u$ as $s\to0$. 
 
 \vs
 
   A partial affirmative answer for the fractional Laplacian is given in \cite{CW0} that for $u \in C^2_c(\R^N)$ and  $x \in \R^N$,
 \begin{equation} \label{eq:limit-behaviour1}
(-\Delta)^\s  u(x) = u(x) + s \cL_1 u (x) + o(\s) \quad{\rm  as}\ \,  \s\to0^+, 
 \end{equation}
where, formally, the operator $\cL_1$ is given as a {\em logarithmic Laplacian}.  Qualitative properties are derived as following.

\begin{proposition}\label{teo-representation}\cite[Theorem 1.1]{CW0}
Let $u \in C^{0,\beta}_c(\R^N)$ for some $\beta\in(0,1)$. Then  for $x \in \R^N$
\begin{eqnarray}
\cL_1  u(x)&:=& \frac{d}{d\s}\Big|_{\s=0} [(-\Delta)^\s u](x)\nonumber
\\[1mm]  &=& c_{N} \int_{\R^N  } \frac{ u(x)1_{B_1(x)}(y)-u(y)}{|x-y|^{N} } dy + \rho_N u(x),   \label{log 1o}  
\end{eqnarray}
where
\begin{equation}
  \label{eq:def-c-N}
c_N:= \pi^{-  N/2}  \Gamma( N/2), \qquad \rho_N:=2 \ln 2 + \psi(\frac{N}{2}) -\gamma
\end{equation}
and $\gamma= -\Gamma'(1)$ is the Euler Mascheroni constant. Here $\psi = \frac{\Gamma'}{\Gamma}$ is the Digamma function. Moreover,
\begin{enumerate}
\item[(i)] for $1 < p \le \infty$, we have $\cL_1  u \in L^p(\R^N)$ and 
\begin{equation}\label{1.1-1}
\frac{(-\Delta)^\s u- u}{\s} \to \cL_1  u \ \, {\rm  in}\ \,  L^p(\R^N)\ \, {\rm as}\ \,  \s \to 0^+ ;
\end{equation}
\item[(ii)] $\mathcal{F}(\cL_1 u)(\xi) = (2 \ln |\xi|)\, \hat u (\xi)$  
 \, for a.e. $\xi \in \R^N$.
\end{enumerate}
\end{proposition}

Following the derivation of the logarithmic Laplacian in \cite{CW0}, it has garnered significant attention as a distinctive nonlocal operator in the investigation of nonlocal problems, including Dirichlet problems \cite{JSW,HS}, related geometric problems   \cite{FKT,CDN}, and Dirichlet eigenvalues \cite{LW,CV,FJW}. Additionally,   the regional integro-differential operator with weak singular kernel could see \cite{TW}.

Motivated by (\ref{1.1-series}), our purpose in this article is to find out  the high order Logarithmic Laplacians  that satisfy the expansions 
$$
\big(\Phi_\s\ast  u\big)(x) =  u(x) + \sum^n_{m=1} \frac{(-s)^m}{m!} \cL_m  u (x) + o(\s^n) \quad  {\rm as}\ \,  \s\to 0^+
$$
and
$$
 (-\Delta )^s u(x) =  u(x) + \sum^n_{m=1} \frac{s^m}{m!} \cL_m  u (x) + o(\s^n) \quad  {\rm as}\ \,  \s\to 0^+.
$$

 To this end,  we first introduce 
some important notations:
 \begin{equation}\label{1.1-k0}
\kappa_1(s) =  \frac{\kappa_{N,s}}{s}= 2^{-2s}\pi^{-\frac N2}  \frac{\Gamma(\frac{N-2s}{2})}{\Gamma(1+s)} 
\end{equation}
and
 \begin{equation}\label{1.1-k1}
\kappa_2(s) = \kappa_{N,s}   \frac{\omega_{_N}}{2s} = 2^{-2s}  \frac{\Gamma(\frac{N-2s}{2})}{\Gamma(\frac N2)\, \Gamma(1+s)},
\end{equation}
where 
$\omega_{_N}$ is the volume of the unit sphere in $\R^N$ with the form
 $\omega_{_N}=\frac{2\pi^{\frac{N}{2}}}{\Gamma(\frac N2)}. $ 
I remark that $\kappa_i\in C^\infty$ in $(-\min\{\frac{N}{2},1\}, \min\{\frac{N}{2},1\})$ with $i=1,2$.
Moreover, we use the notation $\kappa_2^{(m)}(0)$  the $m$-th order  derivative. Compare with (\ref{eq:def-c-N}), we get the relationships
$$\kappa_1(0) =\pi^{-\frac N2}\Gamma(\frac N2)=c_N, \ \ \kappa_2(0)=1 \quad{\rm and}\ \ \kappa_2'(0)=-\rho_N.$$  
For given integer  $n\geq1$,  we denote the kernels
\begin{equation}
  \label{eq:def-jk}
 \k_n(z)=  1_{B_1}(z){\bf q}_n(|z|)\quad\ {\rm and}\quad\ 
  \j_n(z)=   1_{\R^N \setminus B_1}(z){\bf q}_n(|z|), 
\end{equation}
where 
\begin{equation}
  \label{eq:def-k}
  {\bf q}_n(t)=t^{-N}(-2\ln t)^{n-1}\quad{\rm for}\ \, t>0.
   \end{equation}
Given $u\in C^1_c(\R^N)$,  we let $\cK_0  u  =u$ and
\begin{align}
  \cK_n  u (x)  & =   \int_{\R^N}(u(x)-u(y))\k_n(x-y)\,dy - [\j_n \ast u](x) \nonumber
  \\[1mm]& =\int_{\R^N}\big(u(x)1_{B_1(x)}(y)-u(y)\big){\bf q}_n(x-y)\,dy .\label{representation-main}
\end{align}

Our first result show the expansion of $s$-order  Riesz operator $\Phi_{\s} \ast u$  by the $m$-order logarithmic Laplacian  $\cL_m$.  In the following of this paper, we always assume that 
$${\rm  the\ integer} \ \  m\geq1. $$
 Denote by  $C^\alpha_{\rm uloc}(\R^N)$ with $\alpha\in(0,1)$ the set of functions in  $\zeta\in C^{0,\alpha}_{\rm loc}(\R^N)$ such that     
$$\| \zeta \|_{C^\alpha}:=\sup_{x\in\R^N}\| \zeta \|_{C^{0,\alpha}(B_1(x))}<+\infty.  $$

  \begin{theorem}\label{teo 1}
   Let   $\Phi_\s$ be the normalized Riesz kernel  with formula (\ref{1.1-funda}) and 
  $\zeta \in C^\alpha_{\rm uloc}(\R^N)\cap L^\infty(\R^N)\cap L^1(\R^N)$ for some $\alpha \in(0,1)$. 
  For   $x \in \R^N$,  we denote 
\begin{equation} \label{m 1.1}
\cL_m \zeta(x)=    \sum^{m}_{j=0}  \alpha_j \cK_{j} \zeta(x), 
\end{equation} 
where $$\alpha_0=(-1)^{m}  \kappa_2^{(m)}(0),\qquad  \alpha_j=m (-1)^{m+j} \binom{m-1}{j-1} \kappa_1^{(m-j)}(0),\ \ j=1,2,\cdots,m. $$

  Then  we have that 
$$
\frac{d^m}{d\s^m}\Big|_{\s=0^+}  \big(\Phi_{\s} \ast u\big) (x)=  (-1)^m \cL_m  u(x),
$$
\begin{equation} \label{f 1.1}
 \cF\big( \cL_m  u\big)(\zeta)= \big(2\ln |\zeta|\big)^m \, \hat{u}(\zeta)
\end{equation} 
 and 
  for $p\in(1,+\infty]$, $u\in C^{0,\alpha}_c(\R^N)$
$$
\lim_{s\to0^+}\frac{ m!}{\s^{m}} \Big(\Phi_s\ast u -u  -\sum^{m-1}_{j=1}\frac{(-1)^js^j}{j!}   \cL_ju \Big)=(-)^m \cL_mu   \quad {\rm in}\ L^p(\R^N). $$
 
  \end{theorem}

  Here   the binomial coefficient  
    $$\binom{k}{i}=\frac{k!}{i!(k-i)!},$$
where  $i!$ is the $i$  factorial and $0!=1$. 
  \begin{remark}
\begin{enumerate}
\item[(i)] It should be noted that $\cL_m$ is an integro-differential operator with the form (\ref{m 1.1}),  which can be referred to as  the $m$-order logarithmic Laplacian ( $m$ log Laplacian for short)  due to  its Fourier transform (\ref{f 1.1}).   The operators $\cK_j$ are also known as the zero order kernel operators \cite{CS} and they  are related to geometric stable L\'evy processes, as mentioned in reference \cite{KMg,Bass}. 
  
\item[(ii)] $\alpha_m=m\kappa_1 (0)>0$,  so $\cK_m$ is the dominated part of $\cL_m$.  
\item[(iii)] The $m$ logarithmic Laplacian has very weak differential property,   thanks to the singular kernel $\frac{(\ln|x-y|)^{m-1}}{|x-y|^N}$ of $\cK_m$.

\end{enumerate}
\end{remark}

Thanks to the fact that $\Phi_s$ is the fundamental solution of $(-\Delta)^s$, we have the following corollary. 
\begin{corollary}
Let $f\in C^{0,\alpha}_c(\R^N)$, then
the solution $u_{s,f}=\Phi_s\ast f$ of   
\begin{equation}\label{eq 1.1-s poisson}
\left\{ \arraycolsep=1pt
\begin{array}{lll}
\ \ (-\Delta)^s  u=f\quad    {\rm in}\ \,   \R^N,\\[2mm]
 \phantom{    }
 \displaystyle \lim_{|x|\to+\infty} u(x)=0
\end{array}
\right.
\end{equation}
has the expansion that 
$$u_{s,f}=f+\sum_{m=1}^n\frac{(-s)^m}{m!}\cL_mf +o(s^n) \quad {\rm in\ } L^p(\R^N) \ \ {\rm as}\ \ s\to0^+$$
for any $p\in(1,+\infty]$ and any positive integer $n$.

\end{corollary}

 From Theorem \ref{teo 1},  it follows that the $m$ logarithmic Laplacian  is derived from the $m$-order derivative of the order  of Riesz's convolutions.  In fact, it is also  the $m$-order derivative of the order of the fractional Laplacians.

  \begin{theorem}\label{teo 1-flap}
   Let $(-\Delta)^\s$ be the fractional Laplacian with order $\s\in(0,\frac 14)$, 
    $\cL_m$ be the $m$-order Logarithmic Laplacian defined in (\ref{m 1.1}).
   
  Then    for $m\in\N$, $\zeta \in C^{0,\alpha}_c(\R^N)$ with $\alpha>0$ and $x \in \R^N$,
$$
\frac{d^m}{d\s^m}\Big|_{\s=0^+}   (-\Delta)^s \zeta  (x)=    \cL_m  \zeta(x)
$$
and for $p\in(1,+\infty]$ 
$$
\lim_{s\to0^+}\frac{m!}{\s^{m}} \Big((-\Delta)^s\zeta -\zeta -\sum^{m-1}_{j=1}\frac{s^j}{j!}   \cL_j\zeta \Big)=\cL_m\zeta  \quad {\rm in}\ L^p(\R^N). $$
    \end{theorem}

From the above two sides of $m$ derivatives,  we try to enlarge the expansion  (\ref{eq:limit-behaviour1})   to  $m$-order from two sides  of $s\to0$ for any $m\in\N$.
    To this end, we denote   
\begin{equation}\label{1.1-fl-1or}
\cB^{\s}  u(x)  =\left \{
  \begin{array}{lll}
  \big(\Phi_{-\s} \ast u\big)(x)  &\quad  {\rm if} \ \,s\in (-\frac{N}2,0),\\[2mm]
 (-\Delta)^s u(x)  &\quad {\rm if} \ \,s\in (0,+\infty),
  \end{array}
\right.
\end{equation}
where for $s>1$,  $$(-\Delta)^s=(-\Delta)^{s-[\s]}\circ (-\Delta)^{[\s]}  $$
with $[\s]=\sup\{n\in\Z: n\leq \s\}$.    When $s>1$, the higher order Laplacian could also be defined by Fourier transform  
  $$\cF \big((-\Delta)^s u\big)(\zeta)=|\zeta|^{2s}\hat{u}. $$
Note that it is well-defined for $\s\in(-\frac N2,+\infty)$ and    
\begin{equation}\label{1.2-fl-1or}
\cF\big( \cB^{\s}  u\big)(\zeta)=|\zeta|^{2\s}\hat{u}. 
\end{equation}

  The  series states as following.

   \begin{corollary}\label{teo 2}
  Let  $\cL_m$ be the m-order logarithmic operator    given in Theorem  \ref{teo 1}. 
  Then  for  $u \in C^\infty_c(\R^N)$ and  $x \in \R^N$,
$$
\cB^{\s} u(x) =  u(x) + \sum^n_{m=1} \frac{s^m}{m!} \cL_m  u (x) + o(\s^n) \quad  {\rm as}\ \,  \s\to 0.
$$
Moreover, for fixed $s_0\in(-\frac N2,+\infty)$, we have that 
$$
\cB^{\s} u (x) = \cB^{\s_0}  u(x) +\sum^n_{m=1}  \frac{(s-s_0)^{m}}{m!} \cL_m \big(\cB^{\s_0}  u\big) (x) + o((\s-s_0)^n) \quad  {\rm as}\ \,  \s\to s_0.
$$
 \end{corollary}

Note that our expansions in Theorem \ref{teo 2} indicate  that the mapping $\s\to \cB^{\s}u(x)$ is smooth
with respect to the order. 

 Thanks to the singularity of the kernels,   the high order logarithmic Laplacians  have special nonlocal properties and very weak differentiation property.      So there are many new properties for this type of operators related to the order $m$ and our next interest is  to study the fundamental properties: the Dirichlet eigenvalues, basic regularity estimates.  \vs
  
We consider the bounds of the Dirichlet eigenvalues of $\cL_m$
  \begin{equation}\label{eq 1.1}
\left\{ \arraycolsep=1pt
\begin{array}{lll}
\cL_m  u=\lambda  u\quad \  &{\rm in}\ \,   \Omega,\\[1.5mm]
 \phantom{   \cL_m    }
  u=0    &{\rm in}\ \, \R^N\setminus \Omega, 
\end{array}
\right.
\end{equation}
where $\Omega$ is a bounded Lipschitz domain in $\R^N$. 

In order to find these eigenvalues, we need involve the related Hilbert space $\cH_{m,0}(\Omega)$. 
Denote by $C^\infty_c(\Omega)$   the set of functions in $C^\infty(\R^N)$ with compact support in $\Omega\subset \R^N$  and by $\cH_{m,0}(\Omega)$ the closure of functions in $C^\infty_c(\Omega)$  under the norm
$$ \| u \|_{m}=\sqrt{\int_\Omega u^2(x)dx+ \int_{\R^N}\int_{\R^N} (u(x)-u(y))^2 \k_m(x-y) dxdy }. $$
Note that $\cH_{m,0}(\Omega)$ is a Hilbert space with the inner product
$$(u,v)_{\cH_{m,0}}=(u,v)_0+  (u,v)_m,  $$
where
\begin{align*}
(u,v)_0=\int_\Omega uv dx,\quad\  (u,v)_j=   \int_{\R^N}\int_{\R^N} (u(x)-u(y))(v(x)-v(y)) \k_j(x-y)\, dxdy \quad{\rm for}\ \, u,v\in \cH_{m,0}(\Omega).
 \end{align*}

Denote the energy functional of $\cL_m$
  $$\cI_m: \cH_{m,0}(\Omega)\times \cH_{m,0}(\Omega)\to \R,\qquad  \cI_m(u,v) =  \cE_m(u,v)-  \cA_m(u,v),  $$
  where
\begin{align*} 
 \cE_m(u,v) =   \sum^m_{j=0} \alpha_j (u,v)_j,\qquad  \cA_m(u,v) =   \sum^m_{n=1} \alpha_n \int_{\R^N}\int_{\R^N} u(x)v(y)\,\j_n(x-y)  dxdy \quad{\rm for}\ \, u,v\in \cH_{m,0}(\Omega).
\end{align*}
Note that the quadratic form  $\cI_m$  is bilinear, closed, symmetric and semi-bounded in $\cH_{m,0}(\Omega)$.


A function $u\in \cH_{m,0}(\Omega)$ is called the eigenfunction of \eqref{eq 1.1} corresponding to the eigenvalue  $\lambda$ if
\begin{equation}\label{weak-s}
\cI_m(u,\phi)= \lambda\int_{\Omega}u\phi dx,\quad \forall  \, \phi\in \cH_{m,0}(\Omega).
\end{equation}
We have the following characterization of Dirichlet eigenvalues and corresponding eigenfunctions for the operator  $\cL_m$.
 
 \begin{theorem}\label{teo 1-m}
Assume that  $\Omega\subset \R^N$ is a bounded  Lipschitz  domain.  

\begin{itemize}
\item[(i)] Problem \eqref{eq 1.1}  admits  an eigenvalue $\lambda_{m,1}(\Omega)>0$   characterized  by
\begin{equation}\label{Lambda-1-s}
\lambda_{m,1}(\Omega)=\inf_{\substack{u\in \cH_{m,0}(\Omega)\\ u\neq 0}}\frac{\cI_m(u,u)}{\|u\|^2_{L^2(\Omega)}}= \inf_{u\in \cP_{m,1}(\Omega)}\cI_m(u,u),
\end{equation}
with $$\cP_{m,1}(\Omega):=\{u\in \cH_{m,0}(\Omega): \|u\|_{L^2(\Omega)}=1\},$$
 and there exists a  function $\phi_{m,1}\in \cP_{m,1}(\Omega)$ achieving the minimum of $\cI_m$, i.e. 
 $
 \lambda_{m,1}(\Omega)= \cI_{m}(\phi_{m,1},\phi_{m,1}).$ 
 
\item[(ii)]  
  If $m$ is even, then $\lambda_{m,1}(\Omega)\geq0$. 

\item[(iii)] Problem \eqref{eq 1.1}  admits a sequence of eigenvalues $\{\lambda_{m,k}(\Omega)\}_{k\in \N}$ with
\[
-\infty< \lambda_{m,1}(\Omega)\leq \lambda_{m,2}(\Omega)\le \cdots\le \lambda_{m,k}(\Omega)\le \lambda_{m,k+1}(\Omega)\leq \cdots
\]
with corresponding eigenfunctions $\phi_{m,k}$, $k\in \N$ and~
$
\displaystyle \lim_{k\to \infty}\lambda_{m,k}(\Omega) = +\infty.
$

Moreover, for any $k=2,3,\cdots$, the  eigenvalue $\lambda_{m,k}(\Omega)$ can be characterized as
\begin{equation}\label{Lambda-k-s}
\lambda_{m,k}(\Omega)= \inf_{u\in \cP_{m,k}(\Omega)}\cI_m(u,u),
\end{equation}
where 
\[
\cP_{m,k}(\Omega):= \Big\{ u\in \cH_{m,0}(\Omega): \int_{\Omega}u\phi_{m,j}dx  =0 \ \text{ for \ } j=1,2,\cdots k-1 \ \text{ and } \ \|u\|_{L^2(\Omega)}=1\Big\}.
\]
\item[(iv)] The sequence of eigenfunctions  $\{\phi_{m,k}\}_{k\in\N}$ corresponding to eigenvalues $\{\lambda_{m,k}(\Omega)\}$ form a complete orthonormal basis of $L^2(\Omega)$ and  an orthogonal system of $\cH_{m,0}(\Omega)$.

 \item[(v)]   
Assume more that $\Omega\subset B_{\frac12 r_0}$, where  $r_0\in(0,1]$ is such that 
 the  sum of the kernel 
 \begin{equation}\label{hhh}
  \sum^m_{j=1} \alpha_j  (-2\ln t)^{j-1} >0\quad {\rm for}\ \, 0<t<{r_0}. 
  \end{equation}
Then 
 the first eigenvalue $\lambda_{m,1}(\Omega)$ is simple and $\phi_{m,1}$ is positive a.e. in $\Omega$,  that is, if $u\in  \cH_{m,s}(\Omega)$ satisfies \eqref{eq 1.1}  with $\lambda=\lambda_{m,1}(\Omega)$,  then $u=\alpha\phi_{m,1}$ for some $\alpha\in \R$.   Moreover, 
 $$0<\lambda_{m,1}(\Omega)< \lambda_{m,2}(\Omega).$$

\end{itemize}
 \end{theorem} 

\begin{remark}
It is worth noting that we don't obtain the one-fold of  the first Dirichlet eigenvalue of $\cL_m$  and  the positivity of the corresponding eigenfunction in a general bounded domain when $m\geq2$.   The difficulty arises from  $\cL_m=\sum^m_{j=0}\alpha_j\cK_j $ with $\alpha_m>0$,  without the signs of $\{\alpha_j\}$ for $j\leq m-1$, which leads to the failure of  the normal arguments for the positivity of the first Dirichlet eigenvalue. 

 To this end,  we only consider the eigenvalues for small domain $\Omega\subset B_{r_0}$, where $r_0\in(0,\frac12]$  such that  the the sum of the kernel $ \sum^m_{j=0}\alpha_j {\bf q}_j>0$ in $B_{r_0}$. 
  
We also appendix  the Dirichlet eigenvalues of  the dominated part $\cK_m$ of $\cL_m$. In this case,  we can show that first Dirichlet eigenvalue  is simple and the corresponding eigenfunction is positive. 
  
\end{remark}

Finally, we study the  Faber-Krahn inequality of the first eigenvalue for the $m$-order logarithmic Laplacian.   
Recall that 
 $ \displaystyle\sum^m_{j=1} \alpha_j  (-2\ln t)^{j-1}>0$ for $t\in(0,r_0)$ by (\ref{hhh}). 
Let    
$$r_m=\sup\Big\{ r\in(0,r_0):\,  t\to \sum^m_{j=1} \alpha_j  (-2\ln t)^{j-1} \ \text{  is decreasing in $(0,r)$}\Big\},$$ 
 then $r_m\in(0,r_0]$.    

\begin{theorem} \label{sec:faber-Krahn}
   Let $m\geq 2$,  $\Omega$ be a Lipschitz domain such that  $\overline{\Omega}\subset B_{\frac12r_m}$
   and $\bar r\in(0, \frac12r_m)$ be such that  $|\Omega| =|B_{\bar r}|$.
   
   Then  only the ball $B_{\bar r}$  minimizes $\lambda_{m,1}$, i.e. 
$$\lambda_{m,1}(\Omega)\geq \lambda_{m,1}(B_{\bar r}),$$
where $'='$ holds only  for $\Omega=B_{\bar r}$. 


Moreover, the first eigenfunction $\phi_{m,1}$ for $\Omega=B_{\bar r}$ with $\bar r\leq r_m$ is radially symmetric and decreasing with respective to $|x|$.  

\end{theorem} 

\begin{remark}
The  Faber-Krahn inequality of the first eigenfunction holds in the case of the fractional Laplacian,  see  \cite[Theorem 5]{BL},
but it fails for  the regional fractional laplacian.  In fact, \cite{LW} shows that there exists a nonnegative $u\in C^\infty_c(B_1)$
such that
$$\int_{B_1}\int_{B_1} \frac{(u(x)-u(y))^2}{|x-y|^{N+2s}}dxdy<\int_{B_1}\int_{B_1} \frac{(u^*(x)-u^*(y))^2}{|x-y|^{N+2s}} dxdy. $$
The authors in \cite{JKX} provided a weak version of  Faber-Krahn inequality for the regional fractional Laplacian.\vs 

When $\overline\Omega\subset B_{\frac12r_m}$, we see that   
$$\int_{\R^N}\int_{\R^N}\frac{\big(u (x)-u (y)\big)^2}{|x-y|^{N }}{\bf h}_m(|x-y|) \, dxdy  
  = \int_{ B_{r_m}} \int_{B_{r_m}}\frac{\big(u (x)-u (y)\big)^2}{|x-y|^{N}}{\bf h}_{m}(|x-y|)\, dxdy, $$ 
 but we can obtain the  Faber-Krahn inequality for the first eigenfunction. Our method to derive the  Faber-Krahn inequality  based on the rearrangement.\vs
 
  For $m=1$, the Faber-Krahn inequality  is obtained in \cite{CW0} by the relationship 
that $\lambda_{1,1}=\displaystyle\lim_{s\to0^+}\frac{1}{s}(\lambda_{s,1}-1)$, where $\lambda_{s,1}$ is the first Dirichlet eigenvalue of the fractional Laplacian  and by the  Faber-Krahn inequality of the factional Laplacian.  However, for $m\geq2$, it is still open to give an expansion of $\lambda_{s,1}$. 
 For general bounded Lipschtiz domain, it is open to show the  Faber-Krahn inequality for the first eigenfunction
\end{remark}

    The rest of this paper is organized as follows.  Section 2 is devoted to the analysis of the first order of expansion. In Section 3, we derive the second order logarithmic Laplacian and higher order    logarithmic Laplacians are derived by the Riesz operators.   In Section 4,  we obtain  the expansion of   the fractional Laplacians.   Finally, we study the qualitative properties of the operator $\cL_m$ in Section 5 and the Dirichlet eigenvalues for $\cK_m$ are added in the appendix.

 \setcounter{equation}{0}
\section{   Expansion for Riesz convolution }

We fist introduce some notations.  Let $B_r(x)$ be the ball centered at $x\in\R^N$ with radius $r>0$; $B_r=B_r(0)$, 
 $${\bf I}[u ](x,y)=u (y)-u (x)1_{B_1(x)}(y). $$

Observe that 
\begin{equation}
  \label{representation-main*}
  \cK_n  u (x) =   \int_{\R^N} {\bf q}_n(x-y){\bf I}[u ](|x-y|)\,dy  
\end{equation}
 For $s\in(0,\frac N2)$  direct computation shows that for   $x \in \R^N$
\begin{align}\label{2.0-00}
(\Phi_\s\ast \zeta )(x)-\zeta (x)
=\kappa_1(s)s \int_{\R^N} |x-y|^{2s-N}   {\bf I}[\zeta ](x,y) dy +(k_2(s) -1)\zeta (x)
  \end{align}
by  the fact that 
 \begin{align*} 
\kappa_{N,s}  \int_{ B_1(x)} |x-y|^{2s-N} dy =    \kappa_{N,s} \omega_{_N}(2s)^{-1}=2^{-2s}  \frac{\Gamma(\frac{N-2s}{2})}{\Gamma(\frac N2)\, \Gamma(1+s)}=k_2(s).  
  \end{align*}
  Moreover, we see that 
  \begin{eqnarray}\label{e 2.0-1or}  
 \kappa_1(\s)= \frac{\kappa_{N,s}}{s}=2^{-2s}\pi^{-\frac N2}  \frac{\Gamma(\frac{N-2s}{2})}{\Gamma(1+s)}\to \pi^{-\frac N2} \Gamma(\frac{N}{2})= \kappa_1(0)\quad\ {\rm as}\ \, \s\to0^+. 
  \end{eqnarray}

   In this section,  we show that for $\zeta \in C^{0,\alpha}_c(\R^N)$ and  $x \in \R^N$,
$$
\big(\Phi_\s\ast \zeta\big)(x) =  \zeta  +\sum^{m-1}_{j=1}(-1)^j  \frac{s^j}{j!} \cL_j\zeta + o(\s) \quad  {\rm as}\ \,  \s\to 0^+,
$$
where   $\alpha\in(0,1)$ and  $\cL_j \zeta:= (-1)^j \frac{d^j}{d\s^j}\Big|_{\s=0} \Phi_s\ast \zeta$.

  \begin{proposition}\label{teo 1-1or}
   Let $\s\in(0,\frac N2)$ and  $\Phi_\s$ be defined in (\ref{1.1-funda}) and $\cL_m $ be given in (\ref{m 1.1}) for $m\geq1$.\ms
   
  Then 
  $(i)$ for $\zeta \in C^\alpha_{\rm uloc}(\R^N)\cap L^\infty(\R^N)\cap L^1(\R^N)$ with  $\alpha \in(0,1)$,
  $$\lim_{s\to0^+} \frac{m!}{\s^{m}} \Big((\Phi_\s\ast \zeta ) -\zeta  -\sum^{m-1}_{j=1}(-1)^j  \frac{s^j}{j!} \cL_j\zeta   \Big) =(-1)^m  \cL_m\zeta \quad{\rm}\ \, {\rm in}\ L^\infty(\R^N); $$
$(ii)$  for any $p>1$, $\zeta\in C^{0,\alpha}_c(\R^N)$ with  $\alpha \in(0,1)$,
  $$\lim_{s\to0^+} \frac{m!}{\s^{m}} \Big((\Phi_\s\ast \zeta ) -\zeta  -\sum^{m-1}_{j=1}(-1)^j  \frac{s^j}{j!} \cL_j\zeta   \Big) =(-1)^m  \cL_m\zeta \quad{\rm}\ \, {\rm in}\ L^p(\R^N);  $$
  $(iii)$  for $\zeta \in C^{0,\alpha}_c(\R^N)$ with  $\alpha \in(0,1)$,
$$
\frac{d^m}{d\s^m}\Big|_{\s=0^+}  \big(\Phi_{\s} \ast u\big) (x)=  (-1)^m \cL_m u(x)
$$
and
$$
  \cF(\cL_{m} \zeta  )(z) =  (2 \ln |z| )^{m}\,\hat{\zeta} (z) \qquad  {\rm for\ almost\ every}\ z \in \R^N.
$$
  
  \end{proposition}

   \subsection{ Convergence for $m=1$ }

We start the proof by the limit of $ \frac{(\Phi_\s\ast \xi)(x)-\xi(x)}{\s} $ as $\s\to0^+$.  
\begin{lemma}\label{lm 2.1-1or}
Let $\zeta \in C^\alpha_{\rm uloc}(\R^N)\cap L^\infty(\R^N)\cap L^1(\R^N)$ for some $\alpha \in(0,1)$

Then we have that
\begin{eqnarray}\label{e 2}
 \frac{(\Phi_\s\ast \zeta )(x)-\zeta (x)}{\s} \to-\cL_1 \zeta (x) \quad {\rm in}\ L^\infty(\R^N)\ \ \ {\rm  as}\ \, \s\to0^+. 
\end{eqnarray}

\end{lemma}
\noindent{\bf Proof. } For $s\in(0,\min\{\frac14,\frac{N}{2}\})$, direct computation shows that for $\zeta  \in C^\alpha_{\rm uloc} (\R^N)\cap L^1(\R^N, \frac{dx}{(1+|x|)^{N-\alpha }})$ and $x \in \R^N$
\begin{align*}
 \cL_1  \zeta (x) =   \kappa_1(0) \int_{\R^N}  \frac{ {\bf I}[\zeta ](x,y)}{|x-y|^{N} }  dy  -\kappa_2'(0)  \zeta (x).
\end{align*}

\begin{align*}
  \frac{(\Phi_\s\ast \zeta )(x)-\zeta (x)+s\cL_1\zeta (x)}{  \s}  & =  \kappa_1(0)  \int_{\R^N} \Big(\frac1{|x-y|^{N-2s}}-\frac1{|x-y|^{N}}\Big) {\bf I}[\zeta ](x,y)dy  
  \\[1mm]&\quad +  \Big(\kappa_1(\s)-\kappa_1(0) \Big)\int_{\R^N}|x-y|^{2s-N}{\bf I}[\zeta ](x,y)dy  
    \\[1mm]&\quad + \frac{\kappa_2(s)-1-s\kappa_2'(0)} {s} \zeta (x).
\end{align*}

 For $\zeta  \in C^\alpha_{uloc} (\R^N)$, there holds
  $\big|\zeta (y)-\zeta (x)\big|\leq \| \zeta \|_{C^\alpha} |x-y|^\alpha $ for $|x-y|<1$
 and
 \begin{align*} 
 \Big| \int_{B_1(x)}(|x-y|^{2s-N}-|x-y|^{-N}) \big(\zeta (y)-\zeta (x)\big)dy \Big| \leq  \|  \zeta \|_{C^\alpha} \int_{B_1(0)}\frac{ 1-|y|^{2s} }{|y|^{N-\alpha}}  dy.
 \end{align*}
 
Given $\epsilon\in(0,\frac1e)$, we take 
$$\s =  \frac12\frac{\ln(1-\epsilon^{N+2})}{\ln \epsilon} >0,$$ 
where $ \s\sim \frac12\frac{  \epsilon^{N+2}}{-\ln \epsilon} \to 0^+$ as $\epsilon\to0^+$. 
Moreover, we see that 
 $ 1-\epsilon^{2s}= \epsilon^{N+2}$ with the setting of $s$.  

Observe that   
 \begin{align*} 
 \int_{B_\epsilon }\frac{ 1-|y|^{2s} }{|y|^{N-\alpha}} dy\leq  \int_{B_\epsilon } |y|^{\alpha-N}  dy\leq \frac{  \omega_{_N}}{\alpha} \epsilon^\alpha. 
  \end{align*}
 Since
 $$0<1-|y|^{2s}<1-\epsilon^{2s}= \epsilon^{N+2}\quad{\rm  for}\ \, |y|\in(\epsilon,1),$$
then we obtain that 
  \begin{align*} 
   \int_{B_1 \setminus B_\epsilon } \frac{1-|y|^{2s}  }{|y|^{N-\alpha}} dy  \leq \epsilon \int_{B_1 \setminus B_\epsilon } \frac1{|y|^{N-\alpha}} dy\leq  \frac{ \omega_{_N}}{\alpha}   \epsilon^{N+2},
  \end{align*}
  where we recall that $\omega_{_N}$ is the volume of the unit sphere in $\R^N$.
  
  Thus, we have that
   \begin{eqnarray}\label{e 2.0-a}
 \lim_{\s\to0^+} \sup_{x\in\R^N} \int_{B_1(x)} \Big(\frac1{|x-y|^{N-2s}}-\frac1{|x-y|^{N}}\Big) \big(\zeta (y)-\zeta (x)\big)dy =0.
    \end{eqnarray}

Note that $\kappa_1(\s)$ is smooth at $\s=0$, so  there exists  $\s_1>0$ such that  
 \begin{eqnarray}\label{e 2.0-1or-1}  
 |\kappa_1(\s) -\kappa_1(0)|\leq 2|\kappa_1'(0)| \s \quad\ {\rm for}\ \, |\s|\leq \s_1,
  \end{eqnarray}
then 
 \begin{align}
\Big| \Big(\kappa_1(\s)-\kappa_1(0) \Big)\int_{B_1(x)}|x-y|^{2s-N} \big(\zeta (y)-\zeta (x)\big)dy  \Big| &\leq    2|\kappa_1'(0)| \s \, \|  \zeta \|_{C^\alpha} \int_{ B_1(x)} |x-y|^{2s+\alpha -N}  dy 
\nonumber \\ &\leq    2|\kappa_1'(0)| \|  \zeta \|_{C^\alpha} \frac{\omega_{_N}}{\alpha} \s \label{e 2.1-a}
 \\  & \to 0\quad{\rm as}\ \s\to0^+.\nonumber
 \end{align}

 Take  $R=\frac1\epsilon>1$, then we see that 
 \begin{align*}
\Big| \int_{\R^N\setminus B_R(x)}\big(|x-y|^{2s-N}-|x-y|^{-N}\big)  \zeta (y)  dy \Big| &\leq    \int_{\R^N\setminus B_R(x)} |x-y|^{2s-N}  |\zeta (y) | dy 
 \\ &\leq  R^{2s-N} \int_{\R^N\setminus B_R(x)}   |\zeta (y)| dy
 \\ & <   \|\zeta \|_{L^1(\R^N)} \epsilon^{N-\frac12} 
 \end{align*}
and  
 \begin{align*} 
 \Big| \int_{B_R(x)\setminus B_1(x)}\big(|x-y|^{2s-N}-|x-y|^{-N}\big)  \zeta (y)  dy \Big|  &\leq    \int_{B_R(x)\setminus B_1(x)}\big(|x-y|^{2s}-1\big)|x-y|^{-N}  |\zeta (y) | dy 
 \\ &\leq  \|\zeta \|_{L^\infty(\R^N)}\int_{B_R(x)\setminus B_1(x)}\big(|x-y|^{2s}-1\big)      dy 
  \\[1mm]  &\leq \|\zeta \|_{L^\infty(\R^N)} \omega_{_N}  R^{N} (R^{2s}-1)
 \\[1mm]  &< \|\zeta \|_{L^\infty(\R^N)} \omega_{_N}  \epsilon^{-N-2s}(1-\epsilon^{2s})
\\[1mm] &\leq  \|\zeta \|_{L^\infty(\R^N)} \omega_{_N}  \epsilon.
 \end{align*}
 Therefore, 
 \begin{eqnarray}\label{e 2.2-a} 
  \lim_{\s\to0^+} \sup_{x\in\R^N}\Big( \int_{\R^N\setminus B_1(x)}\big(|x-y|^{2s-N}-|x-y|^{-N}\big)  \zeta (y)  dy \Big)=0.
 \end{eqnarray} 
Now we can derive that 
 \begin{eqnarray}\label{e 2.2-1or} 
  \kappa_1(\s) \int_{\R^N\setminus B_1(x)} |x-y|^{2s-N} \zeta (y)dy
\to  \kappa_1(0) \int_{\R^N\setminus B_1(x)} |x-y|^{-N} \zeta (y)dy\quad {\rm in}\ L^\infty(\R^N)
\quad {\rm as}\ \, \s\to0^+ . 
 \end{eqnarray} 
 
 Finally, direct computation shows that
\begin{align*} 
 \kappa_{N,s}  \int_{ B_1(x)} |x-y|^{2s-N} dy-1 &=   \kappa_{N,s} \omega_{_N} (2s)^{-1}-1 
\\[1mm]& =   2^{-2s}\frac{\Gamma(\frac{N-2s}{2})}{\Gamma(\frac N2)\Gamma(1+s)}  -1
  = \kappa_2(s)-1 \to 0 \quad {\rm as}\ \, \s\to0^+ . 
 \end{align*}
which deduces by L'Hospital's rule  that 
  \begin{align*} 
\lim_{s\to0^+}\frac{\kappa_{N,s}  \int_{ B_1(x)} |x-y|^{2s-N} dy-1} {s}&=  \lim_{s\to0^+}\frac{\kappa_{N,s} \omega_{_N}(2s)^{-1}-1} {s} 
\\&  =  \lim_{s\to0^+}     \frac{\kappa_2(s)  -1 } {s}
\\&  =  \kappa_2'(0), 
 \end{align*}
thus,  we have that 
 \begin{eqnarray}\label{e 2.3-a} 
 \lim_{\s\to0^+} \sup_{x\in\R^N}\Big( \frac{\kappa_{N,s}  \int_{ B_1(x)} |x-y|^{2s-N} dy-1} {s} \zeta (x)-\kappa_2'(0)    \zeta (x)\Big)=0. 
 \end{eqnarray} 
 
Combining (\ref{e 2.0-a}), (\ref{e 2.1-a}),  (\ref{e 2.2-a}) and (\ref{e 2.3-a}), we obtain that 
\begin{align*}
  \frac{(\Phi_\s\ast \zeta )(x)-\zeta (x)}{  \s}  \to -\cL_1 \zeta (x)  \quad {\rm in}\ L^\infty(\R^N)\  {\rm as}\ \, \s\to0^+ . 
\end{align*}
   We complete the proof.\hfill$\Box$ 
 
 \subsection{High order expansion }

 \noindent{\bf Proof of Proposition \ref{teo 1-1or}. }  Recall that 
 \begin{eqnarray*}
 \cL_m  \zeta(x) =  (-1)^{m-1} m \sum^{m-1}_{i=0} \binom{m-1}{i}    \kappa_1^{(i)}(0) \int_{\R^N}  \frac{(2\ln |x-y|\big)^{m-1-i}}{|x-y|^{N} } {\bf I}[\zeta](x,y)dy  +(-1)^m   \kappa_2^{(m)}(0)\, \zeta (x). 
 \end{eqnarray*}
  Observe that for $m\geq 2$, 
 \begin{align*} 
 \Theta_{m,s}[\zeta](x)&:=  \frac{1}{\s^{m}} \Big((\Phi_\s\ast \zeta )(x)-\zeta (x)-\sum^{m}_{j=1}(-1)^j  \frac{s^j}{j!} \cL_j\zeta (x) \Big)
  \\[1mm] & =    \frac{\kappa_1(\s)-\kappa_1(0)}{s^{m-1}} \int_{\R^N}|x-y|^{2s-N} {\bf I}[\zeta ](x,y) dy
  \\[1mm]& \quad  + \kappa_1(0) \frac1{\s^{m-1}} \int_{\R^N}\Big(\frac1{|x-y|^{N-2s}}-\frac1{|x-y|^{N}}\Big)  {\bf I}[\zeta ](x,y) dy
  \\[1mm]&\quad  + \frac{\kappa_2(s)-1} {s^{m}} \zeta (x)
   \\[1mm]&\quad-\frac{1}{\s^{m-1}} \sum^{m}_{l=1}\sum^{l-1}_{i=0} \frac{s^{l}  }{ i!(l-1-i)!} \kappa_1^{(i)}(0)  \int_{\R^N}\frac{\big(-2\ln |x-y|\big)^{l-1-i}}{|x-y|^{N} } \, {\bf I}[\zeta ](x,y) dy
  \\[1mm]&\quad - \frac{1}{\s^{m}} \sum_{j=1}^{m} \frac{s^j}{j! }  \kappa_2^{(j)}(0)\, \zeta (x)
  \allowdisplaybreaks
   \\[1mm] & =   \frac1{s^{m-1}} \Big(\kappa_1(\s)-\sum_{j=0}^{m-1} \frac{s^j}{j!}  \kappa_1^{(j)}(0) \Big) \int_{\R^N}|x-y|^{2s-N} {\bf I}[\zeta ](x,y) dy
   \\[1mm]&\quad  + \frac1{s^{m-1}}\sum_{j=0}^{m-1} \frac{s^j}{j!}  \kappa_1^{(j)}(0) \int_{\R^N}|x-y|^{2s-N} {\bf I}[\zeta ](x,y) dy
   \\[1mm]&\quad  - \frac{1}{\s^{m-1}} \sum^{m}_{l=1}\sum^{l-1}_{i=0} \frac{s^{l}  }{ i!(l-1-i)!} \kappa_1^{(i)}(0)  \int_{\R^N}\frac{\big(-2\ln |x-y|\big)^{l-1-i}}{|x-y|^{N} } \, {\bf I}[\zeta ](x,y) dy
   \\[1mm]&\quad   +\frac{1}{s^{m}}\Big(\kappa_2(s)-1-\sum_{j=1}^m \frac{ s^j}{j!} \kappa_2^{(j)}(0) \Big)\zeta (x).
   \end{align*} 
   Note that    
   $$\sum^{m}_{l=1}\sum^{l-1}_{k=0}a_{l,k}=\sum^{m-1}_{k=0}  \sum^{m}_{l=k+1}a_{l,k},  $$
then   
    \begin{align} &\quad\ \frac{1}{\s^{m-1}} \sum^{m}_{l=1}\sum^{l-1}_{k=0} \frac{s^{l}  }{ k!(l-1-k)!} \kappa_1^{(k)}(0)  \int_{\R^N}\frac{\big(-2\ln |x-y|\big)^{l-1-k}}{|x-y|^{N} } \, {\bf I}[\zeta ](x,y) dy\nonumber
    \\[1mm]&=\frac{1}{\s^{m-1}}\sum^{m-1}_{k=0}  \sum^{m}_{l=k+1} \frac{s^{l}  }{ k!(l-1-k)!} \kappa_1^{(k)}(0)  \int_{\R^N}\frac{\big(-2\ln |x-y|\big)^{l-1-k}}{|x-y|^{N} } \, {\bf I}[\zeta ](x,y) dy\label{trans1}
   \\[1mm]&= \frac{1}{\s^{m-1}}\sum^{m-1}_{k=0}  \sum^{m-1-k}_{n=0} \frac{s^{1+k+n}  }{ k!n!} \kappa_1^{(k)}(0)  \int_{\R^N}\frac{\big(-2\ln |x-y|\big)^{n}}{|x-y|^{N} } \, {\bf I}[\zeta ](x,y) dy, \nonumber
      \end{align} 
where we replace $l$ by $n=l-1-k$ in the last equality.   
As a consequence, we derive that 
    \begin{align*} 
  \Theta_{m,s}[\zeta](x)  & =   \frac1{s^{m-1}} \Big(\kappa_1(\s)-\sum_{j=0}^{m-1} \frac{s^j}{j!}  \kappa_1^{(j)}(0) \Big) \int_{\R^N}|x-y|^{2s-N} {\bf I}[\zeta ](x,y) dy
   \\[1mm]&\quad\!  +  \sum_{j=0}^{m-1}  \kappa_1^{(j)}(0)   \int_{\R^N} \frac1{s^{m-1-j} j!}   \Big( |x-y|^{2s}-\sum_{i=0}^{m-1-j} \frac{s^i}{i!}   (2\ln|x-y|)^i  \Big)|x-y|^{-N} {\bf I}[\zeta ](x,y) dy
   \\[1mm]&\quad   +\frac{1}{s^{m}}\Big(\kappa_2(s)-1-\sum_{j=1}^m \frac{ s^j}{j!} \kappa_2^{(j)}(0) \Big)\zeta (x).
\end{align*} 
Our aim is to show that 
 $$\lim_{s\to0^+}\|  \Theta_{m,s}[\zeta]\|_{L^\infty(\R^N)}=0\quad {\rm and}\quad \lim_{s\to0^+}\|  \Theta_{m,s}[\zeta]\|_{L^2(\R^N)}=0.$$

 {\it Part 1: $\displaystyle\lim_{s\to0^+}\|  \Theta_{m,s}[\zeta]\|_{L^\infty(\R^N)}=0$. }  We take  $0<\epsilon<\frac1e$.    For $\zeta  \in C^\alpha_{uloc} (\R^N)$, there holds
  $$\big|\zeta (y)-\zeta (x)\big|\leq \| \zeta \|_{C^\alpha} |x-y|^\alpha \quad {\rm for} \ \ |x-y|<1.$$ 
For $\epsilon\in(0, )$ 
 \begin{align*} 
 &\quad\Big|  \int_{B_\epsilon(x)} \frac1{s^{m-1-j} j!}   \Big( |x-y|^{2s}-\sum_{i=0}^{m-1-j} \frac{s^i}{i!}   (2\ln|x-y|)^i  \Big)|x-y|^{-N} {\bf I}[\zeta ](x,y) dy  \Big| 
 \\& \leq   c_m  \|  \zeta \|_{C^\alpha} \int_{B_\epsilon (0)}\frac{ \big|\ln|y|\big|^{m}  }{|y|^{N-\alpha}}  dy 
  \\& \leq   c_{\alpha,m}  \|  \zeta \|_{C^\alpha}  \epsilon^\alpha (-\ln\epsilon)^m  .
 \end{align*}
For fixed $\epsilon$,  there exists $s_\epsilon\in(0,\frac14)$ such that  for $s\in (0,s_\epsilon)$   
 $$ \frac1{s^{m-1-j}}  \Big| |x-y|^{2s}-\sum_{i=0}^{m-1-j} \frac{s^i}{i!}   (2\ln|x-y|)^i \Big|<\big|2\ln|x-y|\big|^{m}s  \quad{\rm  for}\ \, |y|\in(\epsilon,\frac1\epsilon),$$
then we obtain that
 \begin{align*} 
 &\quad\Big|  \int_{B_1(x)\setminus B_\epsilon(x)} \frac1{s^{m-1-j} j!}   \Big( |x-y|^{2s}-\sum_{i=0}^{m-1-j} \frac{s^i}{i!}   (2\ln|x-y|)^i  \Big)|x-y|^{-N} {\bf I}[\zeta ](x,y) dy  \Big| 
 \\& \leq  \frac2{ j!} \|  \zeta \|_{C^\alpha} \int_{B_1(0)}\frac{ \big|2\ln|y|\big|^{m}  }{|y|^{N-\alpha}}  dy\,s
 \\& =\frac{2^{m+1}}{ j!} \|  \zeta \|_{C^\alpha} \frac{m!}{\alpha^m}\,\epsilon.
 \end{align*}
 Take  $R=\frac1\epsilon>1$, then  for some $c_m>0$
 \begin{align*}
&\quad\ \Big| \int_{\R^N\setminus B_R(x)}\frac1{s^{m-1-j} j!}   \Big( |x-y|^{2s}-\sum_{i=0}^{m-1-j} \frac{s^i}{i!}   (2\ln|x-y|)^i  \Big)|x-y|^{-N}  \zeta (y)  dy  \Big| 
\\&\leq  c_m  \int_{\R^N\setminus B_R(x)} |x-y|^{2s-N}  |\zeta (y) | dy 
 \\ &\leq  c_m R^{2s-N} \int_{\R^N\setminus B_R(x)}   |\zeta (y)| dy
 \\ & <  c_m \|\zeta \|_{L^1(\R^N)} \epsilon^{N-\frac12} 
 \end{align*}
and  
 \begin{align*} 
& \quad\ \Big| \int_{B_R(x)\setminus B_1(x)}\frac1{s^{m-1-j} j!}   \Big( |x-y|^{2s}-\sum_{i=0}^{m-1-j} \frac{s^i}{i!}   (2\ln|x-y|)^i  \Big)|x-y|^{-N}  \zeta (y)  dy \Big|  
\\&\leq   c_m \int_{B_R(x)\setminus B_1(x)}\big(2\ln|x-y|\big)^{m} |x-y|^{-N}  |\zeta (y) | dy\, s
 \\ &\leq c_m \|\zeta \|_{L^\infty(\R^N)}\int_{B_R(x)\setminus B_1(x)} \big(2\ln|x-y|\big)^{m} |x-y|^{-N}     dy\, s
  \\[1mm]  &\leq c_m\|\zeta \|_{L^\infty(\R^N)} \omega_{_N} (\ln R)^{m} \, s
 \\[1mm]  &< c_m\|\zeta \|_{L^\infty(\R^N)} \omega_{_N}  (\ln \frac1\epsilon)^{m} \, \epsilon.
 \end{align*}

By the  arbitrarily of $\epsilon$,  we have that
   \begin{eqnarray}\label{e 2.0-b}
 \lim_{\s\to0^+} \sup_{x\in\R^N} \Big|\int_{\R^N} \frac1{s^{m-1-j} j!}   \Big( |x-y|^{2s}-\sum_{i=0}^{m-1-j} \frac{s^i}{i!}   (2\ln|x-y|)^i  \Big)|x-y|^{-N} {\bf I}[\zeta ](x,y) dy \Big|=0
    \end{eqnarray}

    Note that $\kappa_1(\s)$ is smooth at $\s=0$, so  there exists  $\s_2>0$ such that  
 \begin{eqnarray}\label{e 2.0-1or-1}  
 \Big| \frac1{s^{m-1}} \Big(\kappa_1(\s)-\sum_{j=0}^{m-1} \frac{s^j}{j!}  \kappa_1^{(j)}(0) \Big)\Big|  \leq 2|\kappa_1^{(m)}(0)|\s \quad\ {\rm for}\ \, |\s|\leq \s_2,
  \end{eqnarray}
then 
 \begin{align*}
&\quad\ \Big| \frac1{s^{m-1}} \Big(\kappa_1(\s)-\sum_{j=0}^{m-1} \frac{s^j}{j!}  \kappa_1^{(j)}(0) \Big) \int_{B_1(x)}|x-y|^{2s-N} \big(\zeta (y)-\zeta (x)\big)dy  \Big| 
\\[1mm]  &\leq    2|\kappa_1^{(m)}(0)| \s \, \|  \zeta \|_{C^\alpha} \int_{ B_1(x)} |x-y|^{2s+\alpha -N}  dy 
 \\[1mm]  &\leq    2|\kappa_1^{(m)}(0)| \|  \zeta \|_{C^\alpha} \frac{\omega_{_N}}{\alpha} \s  
 \\[1mm]   & \to 0\quad{\rm as}\ \s\to0^+ 
 \end{align*}
and 
 \begin{align*}
&\quad\ \Big| \frac1{s^{m-1}} \Big(\kappa_1(\s)-\sum_{j=0}^{m-1} \frac{s^j}{j!}  \kappa_1^{(j)}(0) \Big) \int_{\R^N\setminus B_1(x)}|x-y|^{2s-N}  \zeta (y) \big)dy  \Big| 
\\[1mm]  &\leq    2|\kappa_1^{(m)}(0)| \|\zeta\|_{L^1(\R^N)}\s     
 \\[1mm]  & \to 0\quad{\rm as}\ \s\to0^+. 
 \end{align*}
 As a consequence, we obtain that 
    \begin{eqnarray}\label{e 2.1-b}
 \lim_{\s\to0^+} \sup_{x\in\R^N} \Big|\frac1{s^{m-1}} \Big(\kappa_1(\s)-\sum_{j=0}^{m-1} \frac{s^j}{j!}  \kappa_1^{(j)}(0) \Big) \int_{\R^N\setminus B_1(x)}|x-y|^{2s-N}  \zeta (y) \big)dy  \Big| =0.
    \end{eqnarray}

 Finally,  we observe that  
  \begin{align*} 
\lim_{s\to0^+}\frac{  \kappa_2(s)-1-\sum_{j=1}^m \frac{ s^j}{j!} \kappa_2^{(j)}(0)  }{s^{m+1}} = \kappa_2^{(m+1)}(0), 
 \end{align*}
thus,  we have that 
 \begin{eqnarray}\label{e 2.2-b} 
 \lim_{\s\to0^+} \sup_{x\in\R^N}\Big(  \frac{1}{s^{m}}\Big(\kappa_2(s)-1-\sum_{j=1}^m \frac{ s^j}{j!} \kappa_2^{(j)}(0) \Big)\zeta (x) \Big)=0. 
 \end{eqnarray} 
 
Combining (\ref{e 2.0-b}), (\ref{e 2.1-b}) and (\ref{e 2.2-b}), we obtain that 
\begin{align*}
 \lim_{s\to0^+}\|  \Theta_{m,s}[\zeta]\|_{L^\infty(\R^N)}=0 
\end{align*}
and for any $x\in\R^N$
\begin{align}\label{convergence inf}
 \lim_{s\to0^+}  \frac{m!}{\s^{m}} \Big((\Phi_\s\ast \zeta )(x)-\zeta (x)-\sum^{m-1}_{j=1}(-1)^j  \frac{s^j}{j!} \cL_j\zeta (x)\Big)=(-1)^m  \cL_m\zeta (x) 
\end{align}
    
 {\it Part 2: $\displaystyle\lim_{s\to0^+}\|  \Theta_{m,s}[\zeta]\|_{L^p(\R^N)}=0$ for any $p>1$. } Since $C^{0,\alpha}_c(\R^N)\subset \big(C^\alpha_{\rm uloc}(\R^N)\cap L^\infty(\R^N)\cap L^1(\R^N)\big)$, then 
 (\ref{convergence inf}) holds ture for $\zeta\in C^{0,\alpha}_c(\R^N)$. 
 
 For $\zeta\in C^{0,\alpha}_c(\R^N)$,  let 
 $r_0>1$ such that 
 $${\rm supp}(\zeta ) \subset B_{r_0}. $$
Then for $x\in \R^N\setminus B_{2r_0}$,  we see that 
\begin{align*}
\Theta_{m,s}[\zeta](x)  &=     \frac1{s^{m-1}} \Big(\kappa_1(\s)-\kappa_1(0)-\sum_{j=1}^{m} \frac{s^j}{j!}  \kappa_1^{(j)}(0) \Big)  \int_{B_{r_0}(0)\setminus B_1(x)} |x-y|^{2s-N} \zeta (y)dy
  \\[1mm]&\quad \    + \sum_{l=1}^{m}   \kappa_1^{(l)}(0)    \int_{B_{r_0} } \Big(\frac1{|x-y|^{N-2s}}-\sum_{j=0}^{m-1-l} \frac{s^j}{j!} \frac{(2\ln|x-y|)^j}{|x-y|^{N}}\Big) \frac{1}{l!\s^{m-1-l}}  \zeta (y) dy,
\end{align*}
where
\begin{align*} 
 \Big|\sum_{l=1}^{m}   \kappa_1^{(l)}(0)    \int_{\R^N} \Big(\frac1{|x-y|^{N-2s}}-\sum_{j=0}^{m-1-l} \frac{s^j}{j!} \frac{(2\ln|x-y|)^j}{|x-y|^{N}}\Big) \frac{1}{l!\s^{m-l}}  \zeta (y) dy\Big|
  \leq   C \|\zeta _1\|_{L^\infty} |x| ^{-N} \big(|x|^{2\sigma_0}+  (\ln |x|)^m\big)
\end{align*}
for some $C>0$ independent of $s$. 
 For $x\in \R^N\setminus B_{2r_0}$
 \begin{align*}
 \big| \cL_m  \zeta (x)\big| &\leq     2m! \kappa_1(0)\|\zeta \|_{L^\infty}  \int_{B_{r_0}}  \frac{\big|2\ln |x-y|\big|^{m-1}}{|x-y|^{N} }  dy 
 + m! \sum_{i=1}^{m-1}2|\kappa_1^{(i)}(0)|  \|\zeta \|_{L^\infty}  \int_{B_{r_0}}  \frac{\big|2\ln |x-y|\big|^{m-1-i}}{|x-y|^{N} }  dy
 \\[1mm]&\leq C  \|\zeta \|_{L^\infty} |x|^{-N}\big(\ln(|x|)^{m-1}+1\big).
\end{align*}
Thus, we have that 
\begin{align*}
\Big|\frac{1}{\s^{m+1}} \Big(\Phi_\s\ast \zeta )(x)-\zeta (x)-\sum^{m}_{j=1}\frac{s^j}{j!} (-1)^j \cL_j\zeta (x) \Big)   \Big|    \leq C |x|^{2\sigma_0-N}(\ln|x|)^{m},
\end{align*}
where $C>0$ independent of $s$.  Together with the $L^\infty$ bounds, for any $p>1$, we can choose $\sigma_0\in(0,\frac14)$ such that   $(2\sigma_0-N)p<-N$, 
then
\begin{align*}
&\quad\frac{1}{\s^{m}} \Big(\Phi_\s\ast \zeta )(x)-\zeta (x)-\sum^{m}_{j=1}\frac{s^j}{j!} (-1)^j \cL_j\zeta (x) \Big)
\\&\to   \sum^{m}_{l=0} \frac{\kappa_1^{(l)}(0)}{(m-l)! l!}  \int_{\R^N}  \frac{(2\ln |x-y|\big)^{j-l}}{|x-y|^{N} }  {\bf I}[\zeta ](x,y) dy +\frac1{(m+1)!}\, \kappa_2^{(m+1)}(0)\, \zeta (x) \qquad  {\rm in}\ \, L^p(\R^N)\ \ {\rm  as} \ \s\to0^+.
\end{align*}
 That is, $\displaystyle\lim_{s\to0^+}\big\|  \Theta_{m,s}[\zeta]\big\|_{L^p(\R^N)}=0$ for any $p>1$.  \medskip
  
 {\it Part 3: Fourier transform of $\cL_m$. }   
We now have the continuity of the Fourier transform  of  as a map $L^2(\R^N) \to L^2(\R^N)$ and  for $\zeta  \in C^{0,\alpha}_c(\R^N)$,  the Fourier's transform of $\cL_{m}$
 \begin{align*}
\cF(\cL_{m} \zeta ) &=  \lim_{\s \to 0^+} \frac{\cF(\Phi_\s\ast \zeta )(z)-\cF(\zeta )(z)+\sum^{m-1}_{j=1}(-1)^{j}\frac{s^j}{j!} \cF(\cL_j\zeta )(z)}{\frac1{m!}\s^{m}}  
\\[1mm]&=  \lim_{\s \to 0^+}\Bigl(\frac{|z|^{-2\s}-1+ \sum^{m-1}_{j=1}\frac{s^j}{j!}(-1)^{j} (\ln |z|)^j}{\frac1{(m)!}\s^{m}}\Bigr)\cF(\zeta )
\\[1mm]& = (-1)^{m} (2 \ln|z|)^{m} \cF(\zeta ) \quad   {\rm in} \ \,  L^2(\R^N).
\end{align*}
From this we infer that
$$
\cF(\cL_{m} \zeta  )(z) =(-1)^{m}  (2 \ln |z| )^{m}\,\cF(\zeta ) (z ) \qquad  {\rm for\ almost\ every}\ z \in \R^N
$$
 with the formula
 \begin{align*}
 \cL_{m}\zeta (x) = (-1)^{m-1} m \sum^{m}_{l=0} \binom{m-1}{l} \kappa_1^{(l)}(0)   \int_{\R^N}  \frac{(2\ln |x-y|\big)^{m-1-l}}{|x-y|^{N} }{\bf I}[\zeta] (x,y)  dy
   +(-1)^{m}\kappa_2^{(m)}(0)\, \zeta (x) 
 \end{align*}
 as we claimed. \hfill$\Box$

  \medskip

 \noindent{\bf Proof of Theorem \ref{teo 1}. }
   It follows by Proposition \ref{teo 1-1or} directly.  \hfill$\Box$

   \setcounter{equation}{0}
 \section{ Expansions  }
 \subsection{Expansion at $s=0$ for   fractional Laplacian }
For $s\in(0,s)$,  recall
 $$c_{N,s} =2^{2\s}\pi^{-\frac N2} \frac{\Gamma(\frac{N+2\s}2)}{\Gamma(1-\s)}\s.$$
 Then we note that 
 $$c_{N,s} =\kappa_1(-s) s, \quad\quad    c_{N,0}=\kappa_1(0). $$\medskip

 Let 
   $${\bf J}[\zeta ](x,y)=\zeta (x)1_{B_1(x)}(y)-\zeta (y). $$
  \begin{remark}
 We remark that 
 \begin{align*}
 \cL_2  \zeta (x)&=   2\kappa_1(0) \int_{\R^N} |x-y|^{-N} {\bf J}[\zeta] dy
  + 2\kappa_1'(0) \int_{\R^N}  \frac{ 2\ln|x-y| }{|x-y|^{N}}  {\bf J}[\zeta]  dy
   + \kappa_2''(0)  \zeta (x) 
\end{align*}
 and
   \begin{align*}
 \cL_3  \zeta (x)&=    3\kappa_1''(0)     \int_{\R^N}|x-y|^{-N} {\bf J}[\zeta] dy
    + 6  \kappa_1'(0)   \int_{\R^N}  \frac{ 2\ln|x-y| }{|x-y|^{N}}  {\bf J}[\zeta]  dy
   +  3 \kappa_1(0)   \int_{\R^N} \frac{(2\ln|x-y|)^2}{|x-y|^{N}}  {\bf J}[\zeta]  dy
   \\[1mm]&\quad-\kappa_2^{(3)}(0)  \zeta (x).
\end{align*}
  \end{remark}

  \noindent{\bf Proof of Theorem \ref{teo 1-flap}. } 
 Observe that  
   \begin{align*}
 (-\Delta)^\s \zeta (x)-\zeta (x)   & =   \big( \kappa_1(-\s)-\kappa_1(0) \big) \int_{B_1(x)}\frac{{\bf J}[\zeta ](x,y)}{|x-y|^{N+2s} } dy
  \\[1mm]&\quad\  + \kappa_1(0)  \int_{B_1(x)}  \frac{{\bf J}[\zeta ](x,y)}{|x-y|^{N+2s} } dy
 \\[1mm]&\quad\ + \big(\kappa_2(-s)-1\big)\zeta (x),
\end{align*}
where 
$$\int_{\R^N\setminus B_1(x)} \frac1{|x-y|^{N+2s}} dy=\frac{\omega_{_N}}{2s} $$
and
 $$
c_{N,s} \cdot \frac{\omega_{_N}}{2s}  = 2^{2s}  \frac{\Gamma(\frac{N+2s}{2})}{\Gamma(\frac N2)\, \Gamma(1-s)}=  \kappa_2(-s).  $$
  
Direct computations shows that   
 \begin{align*}
 \Lambda_{m,s}[\zeta](x)&:=  \frac{1}{\s^{m}} \Big( (-\Delta)^\s  \zeta (x)-\zeta (x)-\sum^{m}_{j=1}\frac{s^j}{j!} \cL_j\zeta (x) \Big)
  \\[1mm] & =    \frac{\kappa_1(-\s)-\kappa_1(0)}{s^{m-1}} \int_{\R^N}\frac1{ |x-y|^{ N+2s}} {\bf J}[\zeta ](x,y) dy
   \\[1mm]&\quad   + \kappa_1(0)  \int_{\R^N}\frac1{\s^{m-1}}\Big(\frac1{|x-y|^{N+2s}}-\frac1{|x-y|^{N}}\Big)  {\bf J}[\zeta ](x,y) dy
  \\[1mm]&\quad  + \frac{\kappa_2(s)-1} {s^{m}} \zeta (x)
  \\[1mm]&\quad -\frac{1}{\s^{m-1}} \sum^{m}_{j=1}\sum^{j-1}_{i=0} \frac{(-1)^{j-1}s^{i}  }{ i!(j-1-i)!} \kappa_1^{(i)}(0)  \int_{\R^N}\frac{\big(-2\ln |x-y|\big)^{j-i-1} }{  |x-y|^{N}}   {\bf J}[\zeta ](x,y) dy
   \\[1mm]&\quad - \frac{1}{\s^{m}}  \sum_{j=1}^{m} \frac{(-1)^is^j}{j! }  \kappa_2^{(j)}(0)\, \zeta (x)
    \\[1mm] & =  \frac1{s^{m-1}} \Big(\kappa_1(-\s) -\sum_{j=0}^{m} \frac{(-1)^j s^j}{j!}  \kappa_1^{(j)}(0) \Big) \int_{\R^N}|x-y|^{2s-N} {\bf J}[\zeta ](x,y) dy
     \\[1mm] &\quad   +\frac1{s^{m-1}}\sum_{j=0}^{m} \frac{(-1)^j s^j}{j!}  \kappa_1^{(j)}(0) \Big) \int_{\R^N}|x-y|^{2s-N} {\bf J}[\zeta ](x,y) dy
     \\[1mm]&\quad -\frac{1}{\s^{m-1}} \sum^{m}_{j=1}\sum^{j-1}_{i=0} \frac{(-1)^{j-1}s^{i}  }{ i!(j-1-i)!} \kappa_1^{(i)}(0)  \int_{\R^N}\frac{\big(-2\ln |x-y|\big)^{j-i-1} }{  |x-y|^{N}}   {\bf J}[\zeta ](x,y) dy
    \\[1mm]&\quad +\frac{\kappa_2(-s)-1-\sum_{j=1}^m \frac1{j!}(-1)^j s^j \kappa_2^{(j)}(0)} {s^{m}} \zeta (x).
\end{align*}
By (\ref{trans1}) replacing ${\bf I}[\zeta ]$ by ${\bf J}[\zeta ]$,  one has that 
  \begin{align*}
 \Lambda_{m,s}[\zeta](x) &=      \frac1{s^{m-1}} \Big(\kappa_1(-\s) -\sum_{j=0}^{m} \frac{(-1)^j s^j}{j!}  \kappa_1^{(j)}(0) \Big) \int_{\R^N}|x-y|^{2s-N} {\bf J}[\zeta ](x,y) dy
     \\[1mm] &\quad   + \sum_{j=0}^{m} \frac{(-1)^j \kappa_1^{(j)}(0)  }{j! s^{m-1-j}}  \int_{\R^N}
     \Big(|x-y|^{-2s}- \sum^{j-1}_{i=0} \frac{s^{i}  }{ i!} \frac{\big(-2\ln |x-y|\big)^{i} }{  |x-y|^{N}} \Big)|x-y|^{-N} {\bf J}[\zeta ](x,y) dy
    \\[1mm]&\quad +\frac{\kappa_2(-s)-1-\sum_{j=1}^m \frac1{j!}(-1)^j s^j \kappa_2^{(j)}(0)} {s^{m+1}} \zeta (x).
\end{align*}
{\it Part 4: $\displaystyle\lim_{s\to0^+}\|  \Lambda_{m,s}[\zeta]\|_{L^\infty(\R^N)}=0$. }  We take  $0<\epsilon<\frac1e$.    For $\zeta  \in C^\alpha_{uloc} (\R^N)\cap L^\infty(\R^N)\cap L^1(\R^N)$, there holds
  $$\big|\zeta (y)-\zeta (x)\big|\leq \| \zeta \|_{C^\alpha} |x-y|^\alpha \quad {\rm for} \ \ |x-y|<1.$$ 
For $\epsilon\in(0, \frac1e)$ and any $s\in (0,\frac{\alpha}{4})$ 
 \begin{align*} 
 &\quad\Big|  \int_{B_\epsilon(x)} \frac1{s^{m-1-j} j!}   \Big( |x-y|^{-2s}-\sum_{i=0}^{m-1-j} \frac{ s^i}{i!}   (-2\ln|x-y|)^i  \Big)|x-y|^{-N} {\bf I}[\zeta ](x,y) dy  \Big| 
 \\& \leq   c_m  \|  \zeta \|_{C^\alpha} \int_{B_\epsilon (0)}\frac{|y|^{-2s}+ \big|\ln|y|\big|^{m}  }{|y|^{N-\alpha}}  dy 
  \\& \leq   c_{\alpha,m}  \|  \zeta \|_{C^\alpha} \Big(\frac{\epsilon^{\alpha-2s}}{\alpha-2s} +\epsilon^\alpha (-\ln\epsilon)^m \Big). 
 \end{align*}
For fixed $\epsilon$,  there exists $s_\epsilon\in(0,\frac14)$ such that  for $s\in (0,s_\epsilon)$   
 $$ \frac1{s^{m-1-j}}  \Big| |x-y|^{-2s}-\sum_{i=0}^{m-1-j} \frac{s^i}{i!}   (-2\ln|x-y|)^i \Big|<c_4 \Big(1+\big|2\ln|x-y|\big|^{m}\Big)s  \quad{\rm  for}\ \, |y|\in(\epsilon,\frac1\epsilon),$$
 where $c_4 >0$ independent of $x$ and $s$. 
Then we obtain that
 \begin{align*} 
 &\quad\Big|  \int_{B_1(x)\setminus B_\epsilon(x)} \frac1{s^{m-1-j} j!}   \Big( |x-y|^{-2s}-\sum_{i=0}^{m-1-j} \frac{s^i}{i!}   (-2\ln|x-y|)^i  \Big)|x-y|^{-N} {\bf J}[\zeta ](x,y) dy  \Big| 
 \\& \leq  \frac{2c_4}{ j!} \|  \zeta \|_{C^\alpha} \int_{B_1(0)}\frac{ \big|2\ln|y|\big|^{m}  }{|y|^{N-\alpha}}  dy\,s
 \\& =\frac{2^{m+1}c_4}{ j!} \|  \zeta \|_{C^\alpha} \frac{m!}{\alpha^m}\,\epsilon.
 \end{align*}
 Take  $R=\frac1\epsilon>1$, then  for some $c_m>0$
 \begin{align*}
&\quad\ \Big| \int_{\R^N\setminus B_R(x)}\frac1{s^{m-1-j} j!}   \Big( |x-y|^{-2s}-\sum_{i=0}^{m-1-j} \frac{s^i}{i!}   (-2\ln|x-y|)^i  \Big)|x-y|^{-N}  \zeta (y)  dy  \Big| 
\\&\leq  c_5  \int_{\R^N\setminus B_R(x)} |x-y|^{-N} \big(1+ (2\ln|x-y|)^m\big) |\zeta (y) | dy 
 \\ &\leq  c_m R^{-N}\big(1+ (2\ln R)^{m-1}\big) \int_{\R^N\setminus B_R(x)}   |\zeta (y)| dy
 \\ & <  c_m \|\zeta \|_{L^1(\R^N)} \epsilon^{N} (\ln \frac1\epsilon )^m
 \end{align*}
and  
 \begin{align*} 
& \quad\ \Big| \int_{B_R(x)\setminus B_1(x)}\frac1{s^{m-1-j} j!}   \Big( |x-y|^{-2s}-\sum_{i=0}^{m-1-j} \frac{s^i}{i!}   (-2\ln|x-y|)^i  \Big)|x-y|^{-N}  \zeta (y)  dy \Big|  
\\&\leq   c_m \int_{B_R(x)\setminus B_1(x)}\big(2\ln|x-y|\big)^{m} |x-y|^{-N}  |\zeta (y) | dy\, s
 \\ &\leq c_m \|\zeta \|_{L^\infty(\R^N)}\int_{B_R(x)\setminus B_1(x)} \big(2\ln|x-y|\big)^{m-1} |x-y|^{-N}     dy\, s
  \\[1mm]  &\leq c_m\|\zeta \|_{L^\infty(\R^N)} \omega_{_N} (2\ln R)^{m-1} \, s
 \\[1mm]  &< c_m\|\zeta \|_{L^\infty(\R^N)} \omega_{_N}  (\ln \frac1\epsilon)^{m-1} \, \epsilon.
 \end{align*}
  By the  arbitrarily of $\epsilon$,  we have that
   \begin{eqnarray}\label{e 2.0-c}
 \lim_{\s\to0^+} \sup_{x\in\R^N\setminus\{0\}} \Big|\int_{\R^N} \frac1{s^{m-1-j} j!}   \Big( |x-y|^{-2s}-\sum_{i=0}^{m-1-j} \frac{s^i}{i!}   (-2\ln|x-y|)^i  \Big)|x-y|^{-N} {\bf I}[\zeta ](x,y) dy \Big|=0
    \end{eqnarray}

 Furthermore,    it follows by (\ref{e 2.0-1or-1}) that  
 \begin{align*}
&\quad\ \Big| \frac1{s^{m-1}} \Big(\kappa_1(-\s)-\sum_{j=0}^{m-1} \frac{(-1)^js^j}{j!}  \kappa_1^{(j)}(0) \Big) \int_{B_1(x)}|x-y|^{-2s-N} \big(\zeta (y)-\zeta (x)\big)dy  \Big| 
\\[1mm]  &\leq    2|\kappa_1^{(m)}(0)| \s \, \|  \zeta \|_{C^\alpha} \int_{ B_1(x)} |x-y|^{-2s+\alpha -N}  dy 
 \\[1mm]  &\leq    2|\kappa_1^{(m)}(0)| \|  \zeta \|_{C^\alpha} \frac{\omega_{_N}}{\alpha-2s} \s  
 \\[1mm]   & \to 0\quad{\rm as}\ \s\to0^+ 
 \end{align*}
and 
 \begin{align*}
&\quad\ \Big| \frac1{s^{m-1}} \Big(\kappa_1(-\s)-\sum_{j=0}^{m-1} \frac{(-1)^js^j}{j!}  \kappa_1^{(j)}(0) \Big) \int_{\R^N\setminus B_1(x)}|x-y|^{-2s-N}  \zeta (y) \big)dy  \Big| 
\\[1mm]  &\leq    2|\kappa_1^{(m)}(0)| \|\zeta\|_{L^1(\R^N)}\s     
 \\[1mm]  & \to 0\quad{\rm as}\ \s\to0^+. 
 \end{align*}
    \begin{eqnarray}\label{e 2.1-c}
 \lim_{\s\to0^+} \sup_{x\in\R^N} \Big|\frac1{s^{m-1}} \Big(\kappa_1(-\s)-\sum_{j=0}^{m-1} \frac{(-1)^js^j}{j!}  \kappa_1^{(j)}(0) \Big) \int_{\R^N\setminus B_1(x)}|x-y|^{2s-N}  \zeta (y) \big)dy  \Big| =0.
    \end{eqnarray}
Moreover, we see that   
  \begin{align*} 
\lim_{s\to0^+}\frac{  \kappa_2(-s)-1-\sum_{j=1}^m \frac{(-1)^j s^j}{j!} \kappa_2^{(j)}(0)  }{s^{m+1}} =(-1)^{m+1} \kappa_2^{(m+1)}(0), 
 \end{align*}
thus,  we have that 
 \begin{eqnarray}\label{e 2.2-c} 
 \lim_{\s\to0^+} \sup_{x\in\R^N}\Big[  \frac{1}{s^{m}}\Big(\kappa_2(-s)-1-\sum_{j=1}^m \frac{(-1)^j s^j}{j!} \kappa_2^{(j)}(0) \Big)\zeta (x) \Big]=0. 
 \end{eqnarray}

Combining (\ref{e 2.0-c}), (\ref{e 2.1-c}) and (\ref{e 2.2-c}), we obtain that 
\begin{align*}
 \lim_{s\to0^+}\|  \Lambda_{m,s}[\zeta]\|_{L^\infty(\R^N)}=0,
\end{align*}
and for any $x\in\R^N$
$$
 \lim_{s\to0^+}  \frac{1}{\s^{m}} \Big( (-\Delta)^s\zeta (x)-\zeta (x)-\sum^{m-1}_{j=1}(-1)^j  \frac{s^j}{j!} \cL_j\zeta (x)\Big)= \cL_m\zeta (x) 
$$

 {\it Part 5: $\displaystyle\lim_{s\to0^+}\|  \Lambda_{m,s}[\zeta]\|_{L^p(\R^N)}=0$ for any $p>1$. } Since $C^{0,\alpha}_c(\R^N)\subset \big(C^\alpha_{\rm uloc}(\R^N)\cap L^\infty(\R^N)\cap L^1(\R^N)\big)$, then 
 (\ref{e 2.2-c}) holds ture for $\zeta\in C^{0,\alpha}_c(\R^N)$. 
 
 For $\zeta\in C^{0,\alpha}_c(\R^N)$,  let 
 $r_0>1$ such that 
 $${\rm supp}(\zeta ) \subset B_{r_0}. $$
Then for $x\in \R^N\setminus B_{2r_0}$,  we see that 
\begin{align*}
\Lambda_{m,s}[\zeta](x)  &=     \frac1{s^{m-1}} \Big(\kappa_1(-\s)-\kappa_1(0)-\sum_{j=1}^{m} \frac{(-1)^j s^j}{j!}  \kappa_1^{(j)}(0) \Big)  \int_{B_{r_0}(0)\setminus B_1(x)} |x-y|^{-2s-N} \zeta (y)dy
  \\[1mm]&\quad \    + \sum_{l=1}^{m}   \kappa_1^{(l)}(0)    \int_{B_{r_0} } \Big(\frac1{|x-y|^{N+2s}}-\sum_{j=0}^{m-1-l} \frac{s^j}{j!} \frac{(-2\ln|x-y|)^j}{|x-y|^{N}}\Big) \frac{1}{l!\s^{m-1-l}}  \zeta (y) dy,
\end{align*}
where
\begin{align*} 
 \Big|\sum_{l=1}^{m}   \kappa_1^{(l)}(0)    \int_{\R^N} \Big(\frac1{|x-y|^{N+2s}}-\sum_{j=0}^{m-1-l} \frac{s^j}{j!} \frac{(-2\ln|x-y|)^j}{|x-y|^{N}}\Big) \frac{1}{l!\s^{m-l}}  \zeta (y) dy\Big|
  \leq   C \|\zeta _1\|_{L^\infty} |x| ^{-N} \big(1+  (\ln |x|)^m\big)
\end{align*}
for some $C>0$ independent of $s$. 
 For $x\in \R^N\setminus B_{2r_0}$
 \begin{align*}
 \big| \cL_m  \zeta (x)\big| &\leq     C  \|\zeta \|_{L^\infty} |x|^{-N}\big(\ln(|x|)^{m-1}+1\big).
\end{align*}
Thus, we have that 
\begin{align*}
\Big|\frac{1}{\s^{m+1}} \Big( (-\Delta)^s\zeta(x)-\zeta (x)-\sum^{m}_{j=1}\frac{s^j}{j!} (-1)^j \cL_j\zeta (x) \Big)   \Big|    \leq C |x|^{ -N}(\ln|x|)^{m},
\end{align*}
where $C>0$ independent of $s$.  Together with the $L^\infty$ bounds, for any $p>1$, we can choose $\sigma_0\in(0,\frac14)$ such that   $(2\sigma_0-N)p<-N$, 
then
\begin{align*}
 \quad\frac{1}{\s^{m}} \Big((-\Delta)^s\zeta (x)-\zeta (x)-\sum^{m-1}_{j=1}\frac{s^j}{j!}   \cL_j\zeta (x) \Big)
 \to   \cL_m \zeta (x) \qquad  {\rm in}\ \, L^p(\R^N)\ \ {\rm  as} \ \s\to0^+
\end{align*}
  for any $p>1$.   \hfill$\Box$

  \subsection{Taylor's expansion for fractional Laplacians at $s>0$}

 We remark that, particularly,  the Taylor's expansions: for $u \in C^\infty_c(\R^N)$ and  $x \in \R^N$,
$$
(-\Delta)^\s  u(x) =   u(x) +  s   \cL_1 u (x) + o(\s) \quad  {\rm as}\ \,  \s\to 0,
$$     
$$
(-\Delta)^\s  u(x) =   (-\Delta)^{s_0} u(x) +  (s-s_0) \cL_1 \circ (-\Delta)^{s_0}  u (x) + o(\s-{s_0}) \quad  {\rm as}\ \,  \s\to {s_0}
$$  
 and
 $$
\Phi_\s\ast  u(x) =   \Phi_{\s_0}\ast  u(x) -  (s-s_0) \cL_1 \big(\Phi_{\s_0}  \ast u\big) (x) + o(\s-s_0) \quad  {\rm as}\ \,  \s\to s_0, 
$$  
 where ${s_0}$ is a positive integer and $s_0\in(0,\frac N2)$. \ms

 \noindent{\bf Proof of Corollary \ref{teo 2}. }  
   Given $s_0\in(0,\frac N2)$,   for any $u \in C_c^\infty(\R^N)$, 
 $\Phi_{\s_0}\ast  u \in C^\infty(\R^N)$ satisfies 
  $$\Big|\big(\Phi_{\s_0}\ast  u\big)(x)\Big|\leq (1+|x|)^{2 \s_0 -N}\ \,  {\rm for}\ x\in\R^N. $$

Note that for $\s\in \big( -\frac{N}{2}, \s_0\big)$ 
  $$\cF(\cB^{\s} u)=|\cdot|^{2\s}\hat{u}= |\cdot|^{2(\s-\s_0)} |\cdot|^{2\s_0} \hat{u}= \cF(\Phi_{-(\s-\s_0)})  \cdot \cF(  \Phi_{\s_0} \ast u),$$
 then we have that 
 $$\cB^{\s}  u=\Phi_{\s-\s_0}\ast \big(\cB^{\s_0} \ast u).$$
For $\s\in \big(  \s_0, \s_0+\frac12\big)$, we have that 
    $$\cF(\cB^{\s} u)=|\cdot|^{2\s}\hat{u}= |\cdot|^{2(\s-\s_0)} |\cdot|^{2\s_0} \hat{u}= \cF((-\Delta)^{\s-\s_0})   \cF( \big(\cB^{\s_0}   u). $$

 Now we apply  Theorem \ref{teo 1} and Theorem \ref{teo 1-flap} to obtain that 
 $$
 \cB^s  u (x) = \cB^{\s_0}  u (x) +\sum^n_{m=1} (-1)^m (s-s_0)^m\cL_m  \big(\cB^{\s_0}     u\big) (x) + o((\s-s_0)^n) \quad  {\rm as}\ \,  \s\to s_0.
$$
Particular, for $s_0=0$, we have that 
 \begin{align*} 
\cB^s(x)  = \zeta (x) + \sum^{n}_{m=1} \frac{s^m}{m!}\cL_m \zeta  (x) + o(\s^n) \quad  {\rm as}\ \,  \s\to0.
 \end{align*} 
We complete the proof. \hfill$\Box$

     \setcounter{equation}{0}
  \section{Qualitative properties of $\cL_m$}
  \subsection{$m$-order Dini Continuity}
   
 Since $\cK_m$ is  the main part of $\cL_m$, we concentrate on the Dini continuity for the operator $\cK_m$.  We first analyze for which functions $u$ and points $x \in \R^N$, the expression $\cK_m u(x)$ is well-defined pointwisely by the formula (\ref{representation-main}).  

 Let $\Omega \subset \R^N$ be a measurable subset and $u: \Omega \to \R$ be a measurable function. The module of continuity of $u$ at a point $x \in \Omega$
is defined by
$$
\omega_{u,x,\Omega}: (0,\infty) \to [0,\infty),\qquad  \omega_{u,x,\Omega}(r)= \sup_{\stackrel{
y \in \Omega}{|y-x|\le r}} |u(y)-u(x)|.
$$
The function $u$ is called {\em $m$-order Dini continuous} at $x$ if 
$$\displaystyle \int_0^1 \frac{\omega_{u,x,\Omega}(r)}{r}\big(1-\ln r\big)^{m-1}\,dr < +\infty.$$
If
\begin{equation*}
\int_0^1 \frac{\omega_{u,\Omega}(r)}{r}\big(1-\ln r\big)^{m-1}\,dr < +\infty \quad \text{for the uniform continuity module $\omega_{u,\Omega}(r):=
\sup \limits_{x \in \Omega}\omega_{u,x,\Omega}(r)$,}
\end{equation*}
then we call $u$ {\em uniformly  $m$-order Dini continuous} in $\Omega$.  

For $s, t \in \R$,  we also denote $L^1_{s,t}(\R^N)$  the space of locally integrable functions $u: \R^N \to \R$ such that
$$
\|u\|_{L^1_{s,t}}:= \int_0^\infty \frac{|u(x)|}{(1+|x|)^{N+2s}} \big(\ln (e+|x|)\big)^{t-1}\,ds < +\infty.
$$
 
  \begin{proposition}
\label{dini-continuity-well-defined}
Let the integer  $m\geq 1$ and $u \in L^1_{0,t}(\R^N)$ with $t\geq m$.
\begin{enumerate}
\item[(i)] If $u$ is $m$-order Dini continuous at some $x \in \R^N$, then $[\cK_m u](x)$ is well-defined by the formula~(\ref{representation-main}). Moreover, if $\Omega \subset \R^N$ is an open subset and $x \in \Omega$, then we have the alternative representation
\begin{equation}
 \cK_m  u(x) =  \int_{\Omega} \big( u(x) -u(y)\big) {\bf q}_m(|x-y|)  dy -   \int_{\R^N \setminus \Omega} u(y)  {\bf q}_m(|x-y|) dy + h_{m,\Omega}(x) u(x), \label{representation-regional}
\end{equation}
where 
$$h_{m,\Omega}(x)=\int_{B_1(x)\setminus \Omega}{\bf q}_m(|x-y|)   dy-\int_{\Omega \setminus B_1(x)}{\bf q}_m(|x-y|)  dy.$$  
\item[(ii)] If $u$ is uniformly $n$-order  Dini continuous in  $\Omega$ for $n\geq m$, then $\cK_m u$ is continuous in $\Omega$.
\end{enumerate}
\end{proposition}

  To prove this proposition,  we need the following observation.

\begin{lemma}
\label{continuity-convolution}
Let $m\geq1$,  $u \in L^1_{0,m}(\R^N)$ and  $v: \R^N \to \R$ be measurable with $|v(x)| \le C (1+|x|)^{-N}\big(\ln (e+|x|)\big)^{m-1}$ for $x \in \R^N$ with some $C>0$.
Then the convolution $v\ast  u: \R^N \to \R$ is well-defined and continuous.
\end{lemma}
 {\bf Proof. }
This is a direct consequence of Lebesgue's theorem if
$u \in L^1_{0,m}(\R^N)$. Moreover, if $\cO \subset \R^N$ is compact, it holds that for $x \in \cO$
\begin{align*}
\big|[v \ast  u](x)\big| \le \int_{\R^N}|v(x-y)||u(y)|\,dy &\le C \int_{\R^N} \frac{|u(y)|}{(1+|x-y|)^N}\big(\ln (e+|x-y|)\big)^{m-1}\,dy \\
&\le C_{\cO} \int_{\R^N} \frac{|u(y)|}{(1+|y|)^N}\big(\ln (e+|x|)\big)^{m-1}\,dy 
\\[1mm]& = C_{\cO}  \|u\|_{L^1_{0,m}}
\end{align*}
for a constant $C_{\cO}>0$. 
Therefore $\|v \ast  u\|_{L^\infty(\cO)} \le C_{\cO}  \|u\|_{L^1_{0,m}}$. Hence, using the fact that $C_c(\R^N)$ is dense in $L^1_{0,m}(\R^N)$, a standard approximation argument shows that $v\ast u$ is continuous on $\cO$ for arbitrary $u \in L^1_{0,m}(\R^N)$. Since $\cO$ is chosen arbitrarily, we conclude that $v\ast u: \R^N \to \R$ is continuous. \hfill$\Box$\medskip

\noindent {\bf Proof of Proposition \ref{dini-continuity-well-defined}. } For $m=1$, it could see \cite[Proposition 2.2]{CW0}. In the following proof, we take $t=m\geq 2$. 

$(i)$ Since
\begin{align*}
&\quad \int_{\R^N  }(|u(x)1_{B_1(x)}(y)-u(y)|)  {\bf q}_m(|x-y|)   dy
\\ &= \int_{B_1(x)} (|u(x)1_{B_1(x)}(y)-u(y)|) {\bf q}_m(|x-y|)  dy + \int_{\R^N \setminus B_1(x)}
(|u(x)1_{B_1(x)}(y)-u(y)|) {\bf q}_m(|x-y|)  dy\\
&\le 2^{m-1} \omega_{_N} \int_{0}^1 \frac{\omega_{u,x,\Omega}(r)}{r} |\ln r|^{m-1} dr + 2^{m-1} C_x \int_{\R^N}\frac{|u(y)|}{(1+|y|)^N}\big(\ln(e+ |y|)\big)^{m-1}dy
\\& < +\infty
\end{align*}
with a constant $C_x>0$ by assumption, it follows that $ \cK_m  u(x)$ is well-defined by (\ref{representation-main}).  Next, we let $x\in \Omega \subset \R^N$, then  from (\ref{representation-main})  we see that
\begin{align*}
  \cK_m  u(x) &=  \int_{\R^N  }(|u(x)1_{B_1(x)}(y)-u(y)|) {\bf q}_m(|x-y|)  dy \\
&= \int_{\Omega } (|u(x)1_{B_1(x)}(y)-u(y)|) {\bf q}_m(|x-y|)  dy +  \int_{\R^N \setminus \Omega }(|u(x)1_{B_1(x)}(y)-u(y)|) {\bf q}_m(|x-y|) dy \\
&= \int_{\Omega } \big( u(x)-u(y)\big) {\bf q}_m(|x-y|) dy   - u(x) \int_{\Omega \setminus B_1(x)}  {\bf q}_m(|x-y|) dy
\\&\quad + u(x) \int_{B_1(x) \setminus \Omega} {\bf q}_m(|x-y|) dy-  \int_{\R^N \setminus \Omega}  u(y) {\bf q}_m(|x-y|)dy
\\&= \int_{\Omega } \big(u(x)-u(y)\big)  {\bf q}_m(|x-y|) dy -  \int_{\R^N \setminus \Omega}  u(y) {\bf q}_m(|x-y|) dy +  h_{m,\Omega}(x) u(x),
\end{align*}
which yields (\ref{representation-regional}). \smallskip

$(ii)$  We write $\cK_m u =   f_0 - \,\j_m \ast u(x)$, where $\,\j_m$ is defined in (\ref{eq:def-jk}), 
$$
   f_0(x)= \int_{B_1(x)} \big( u(x) -u(y)\big) {\bf q}_m(|x-y|)  dy. 
$$
By Lemma~\ref{continuity-convolution}, we see that $\,\j_m \ast u$ is continuous on $\Omega$. Now we show the continuity of $f_0$ in $\Omega$. 

Let
$\eps_k:= 2^{-2^k}$ for $k \in \N_0$ and  fix  $x \in \Omega$,
 $r:= \min \{\frac14, \frac{\dist(x,\partial \Omega)}{4}\}$.  
 Now we choose  $k \in \N$   such that $\eps_k <r$.  For $y \in \Omega$ with $|x-y|<\eps_k$, we then have
\begin{align*}
&|f_0(x)-f_0(y)| = \Bigl|\int_{B_1} \Big(u(x)-u(x+z) -[u(y)-u(y+z)]\Big)  {\bf q}_m(|z|) dz \Bigr|\\
&\le \int_{B_{\eps_k}}\Big(|u(x)-u(x+z)|+|u(y)-u(y+z)| \Big)  {\bf q}_m(|z|) dz 
\\&\quad +\Bigl| \int_{B_1 \setminus B_{\eps_k}}\Big(u(x)-u(x+z) -[u(y)-u(y+z)]\Big)  {\bf q}_m(|z|) dz\Bigr|\\
&\le 2 \int_{B_{\eps_k}}  \omega_{u,\Omega}(|z|){\bf q}_m(|z|) dz + \omega_{u,\Omega}(|x-y|) \int_{B_1 \setminus B_{\eps_k}}{\bf q}_m(|z|) dz
\\&\quad + \Bigl| \int_{B_1 \setminus B_{\eps_k}} \Big( u(x+z)-u(y+z)\Big){\bf q}_m(|z|) dz\Bigr|
\\
&\le \omega_{_N}\delta_k +\Bigl|\int_{B_1 \setminus B_{\eps_k}}\big(u(x+z)-u(y+z)\big) {\bf q}_m(|z|) dz\Bigr|, 
\end{align*}
where
 $$\delta_k := 2\int_{0}^{\eps_k}  \omega_{u,\Omega}(\tau) {\bf q}_m(\tau)    \,d\tau  +(m-1)! \, \omega_{u,\Omega}(\eps_k) (\ln \frac{1}{\eps_k})^{m}.$$

Now we claim that  $\delta_k \to 0$ as $k\to+\infty$. In fact, we see that 
\begin{align*}
  \sum_{k=1}^\infty \omega_{u,\Omega}(\eps_k) (\ln \frac{1}{\eps_k})^m
&\leq \sum_{k=1}^\infty \omega_{u,\Omega}(\eps_k) \sum^{m-1}_{i=0}  \frac{(m-1)!}{i!}  (\ln \frac{1}{\eps_k})^{i} \ln  \frac{\eps_{k-1}}{\eps_k}
\\&\le   \sum_{k=1}^\infty \int_{\eps_k}^{\eps_{k-1}}\frac{\omega_{u,\Omega}(r)}{r} (\ln \frac1r)^{m-1}\,dr
\\&= 2 \int_{0}^{\frac{1}{4}}  \omega_{u,\Omega}(r)  {\bf q}_m(r)\,dr 
\\&<+\infty
\end{align*}
by assumption of the uniformly $m$-order  Dini continuity. Hence
\begin{equation}
\label{prelim-remark}
\omega_{u,\Omega}(\eps_k) (\ln \frac{1}{\eps_k})^m \to 0 \quad \text{as $k \to \infty.$}
\end{equation}
Thus,  $\delta_k \to 0$  as $k\to +\infty$. 

Finally, we show  $\Bigl| \int_{B_1 \setminus B_{\eps_k}} \Big( u(x+z)-u(y+z)\Big){\bf q}_m(|z|) dz\Bigr|\to0$ if $|x-y|\to0$. 
Denote
$$
  v_k(z)= 1_{B_1 \setminus B_{\eps_k}}(z){\bf q}_m(|z|),
$$
then we apply by Lemma~\ref{continuity-convolution} to obtain that 
$$
\Bigl|\int_{B_1 \setminus B_{\eps_k}}\big(u(x+z)-u(y+z)\big) {\bf q}_m(|z|) dz\Bigr|= \bigl|[v_k* u](x)-[v_k* u](y)\bigr|
\to 0 \qquad \text{as $y \to x$}
$$
for every $k \in \N$. We thus conclude that
$$
\limsup_{|y-x| \to 0^+}\bigl|f_0(x)-f_0(y)\bigr| \le  \omega_{_N}\delta_k \qquad \text{for every $k \in \N$,}
$$
and this implies that $\lim \limits_{|y-x| \to 0^+}|f_0(x)-f_0(y)|= 0$. Hence $f_0$ is continuous in $x$,  so is $\cK_m u$.
 \hfill$\Box$ 
 \begin{corollary} \label{cr dini-continuity-well-defined}
Let the integer  $m\geq 1$ and $u \in L^1_{0,t}(\R^N)$ with $t\geq m$.
\begin{enumerate}
\item[(i)] If $u$ is $m$-order Dini continuous at some $x \in \R^N$, then $ \cL_m u (x)$ is well-defined. 
\item[(ii)] If $u$ is uniformly $m$-order  Dini continuous in  $\Omega$, then $\cL_m u$ is continuous in $\Omega$.
\end{enumerate}
\end{corollary}
{\bf Proof. }   If $\cK_m u (x)$ is well-defined, so is $\cK_n u (x)$ for all $0\leq n< m$. Moreover, if $\cK_m u$ is continuous in $\Omega$, so is $\cK_n u  $ for all $0\leq n< m$.

Our arguments follow by Proposition \ref{dini-continuity-well-defined} and expression of $\cL_m$ in (\ref{m 1.1}).  \hfill$\Box$

  \subsection{Eigenvalues in bounded domains}
  
We start  this subsection from the estimates of related Hilbert space. 
 
 \begin{lemma} Let $m\geq1$, $u\in  \cH_{m,0}(\Omega) $ and $0\leq j< m$, then  for any $\epsilon\in(0,1]$,  there exists $C_j=C_j(\epsilon)>0$ such that  
  \begin{equation} \label{ineq 1}
  (u,u)_j \leq \epsilon \|u\|_m^2+C_j\|u\|_{L^2(\Omega)}
  \end{equation}

  \end{lemma}
  {\bf Proof. } Let $m\geq 2$ and $j=1,\cdots, m-1$, then for given $\epsilon\in(0,1]$ , we see that 
 $$(-2\ln |z|)^j\leq \epsilon (-2\ln |z|)^m\quad {in}\ B_{r_j},$$
 where $r_j=e^{-\frac12 \epsilon^{j-m}}$.

  Direct computation shows that 
\begin{align*}
  (u,u)_j &\leq  \int_{\R^N} \int_{\R^N} \big( u(x)-u(y)\big)^2 {\bf q}_m(x-y)1_{ B_{r_j}}(x-y) dxdy
  \\& \quad + \int_{\R^N} \int_{\R^N} \big( u(x)-u(y)\big)^2 {\bf q}_j(x-y)1_{B_1\setminus B_{r_j}}(x-y) dxdy
\\&\leq  \epsilon \|u\|_m^2+ 2 \int_{\Omega} u(x)^2 \int_{\R^N} {\bf q}_j(x-y)1_{B_1\setminus B_{r_0}}(x-y) dy-
2\int_{\Omega}\int_{\Omega} u(x)u(y) {\bf q}_j(x-y)1_{B_1\setminus B_{r_0}}(x-y) dxdy
\\&\leq \epsilon \|u\|_m^2+  c_{1,j}  \int_{\Omega} u(x)^2dx+  2 \int_\Omega |u(x)| \Big(\int_{B_1(x)\setminus B_{r_0}(x)} u(y) ^2dy\Big)^\frac12  \Big(\int_{B_1(x)\setminus B_{r_0}(x)} {\bf q}_j(x-y)^2 dy\Big)^\frac12 dx
\\&\leq\epsilon \|u\|_m^2+  c_{1,j}  \int_{\Omega} u(x)^2dx+c_{2,j}(r_0)  \Big(\int_\Omega u(x)^2 dx\Big)^\frac12 \Big( \int_\Omega \int_{B_1(x)\setminus B_{r_0}(x)} u(y)^2 dy  dx\Big)^\frac12 
\\&\leq \epsilon \|u\|_m^2+  \Big( c_{1,j}(r_0)  +c_{2,j}\sqrt{ |\Omega|}\Big)   \int_\Omega u(x)^2 dx
\\&\leq \epsilon \|u\|_m^2+ C_j(\epsilon) \int_\Omega u(x)^2 dx,
\end{align*}
where $C_j\leq  c_{1,j}(r_0)  +c_{2,j}(r_j)\sqrt{|\Omega|}$, 
 $$c_{1,j} =2 \int_{B_1\setminus B_{r_j} } {\bf q}_j(y)  dy  =2^{j}\omega_{_N} \int_{r_j}^1  t^{-1} (-\ln t)^{j-1}dt    $$
 and 
  \begin{align*}
  c_{2,j} =2\Big(\int_{B_1\setminus B_{r_j}} {\bf q}_j(y)^2 dy \Big)^\frac12
  \leq 2^{j}   \omega_{_N}^\frac12 r_j^{-\frac N2} \Big( \int_{r_j}^1t^{-1} (-\ln t)^{2j-2} dt\Big)^\frac12.  
  \end{align*}
 We complete the proof. \hfill$\Box$\medskip
 
From above lemma,  for integer $n<m$, a corollary follows that 
\begin{equation}
  \label{eq: embedding}
 \cH_{m,0}(\Omega) \varsubsetneqq  \cH_{n,0}(\Omega)\quad {\rm for}\ n<m. 
  \end{equation}

\begin{proposition}\label{cr 6.1}
  Let $m\geq 1$ and $\Omega$ be a bounded Lipschitz domain.  
  
  $(i)$ If  $u\in \cH_{m,0}(\Omega)$, then $u_{\pm}, |u| \in \cH_{m,0}(\Omega)$,  
  where $u_{\pm}=\pm \max\{0,\pm u\}$. \smallskip

$(ii)$ If  $u\in \cH_{m,0}(\Omega)$, then $\||u|\|_m\leq \| u\|_m$,
where the equality holds only when $u$ does not change sign. 
\smallskip
  
$(iii)$  The embedding $\cH_{m,0}(\Omega) \hookrightarrow L^2(\Omega)$ is continuous and compact and
  there exists $C_m\in(0,1]$ such that 
\begin{equation}
  \label{eq:embedding-compact}
  \|u\|_{L^2(\Omega)} \leq C_m \| u \|_{m}\quad {\rm for\ any}\  u\in \cH_{m,0}(\Omega). 
\end{equation}

\end{proposition}
{\bf Proof.}    Note that  $\cH_{m,0}(\Omega)$ could be also seen the closure of functions in $C^{0,\frac12}_c(\Omega)$  under the norm  $\| u \|_{m}$.  For any $u\in C^{0,\frac12}_c(\Omega)$,  we see that $u_{\pm}\in C^{0,\frac12}_c(\Omega)$, so is $|u|$,  Thus, $|u| \in \cH_{m,0}(\Omega)$. \smallskip

We see that $(u,u)_0=(|u|, |u|)_0$ and for $j\geq 1$, 
\begin{align*}
(|u|,|u|)_j&=    \int_{\R^N}  \int_{\R^N} (|u(x)|-|u(y)|)^2 \k_j(x-y) dxdy  
\\&\leq \int_{\R^N}  \int_{\R^N} (u(x)-u(y))^2 \k_j(x-y) dxdy=(u,u)_j,
\end{align*}
where the equality holds only when $u$ does not change sign. 
Our argument $(ii)$ follows by the definition of $\| u \|_{m}^2=(u,u)_0+(u,u)_m$. \smallskip

  The kernel ${\bf q}_m$ verifies the assumptions of \cite[Theorem 2.1]{CP}, then it implies  that  
 the embedding $\cH_{m,0}(\Omega) \hookrightarrow L^2(\Omega)$ is continuous and compact 
and (\ref{eq:embedding-compact}) is obvious by the definition of the norm of $\cH_{m,0}$.    \hfill$\Box$\medskip

Let ${\bf H}_{m,0}(\Omega)$ be
 the completion of $C^\infty_c(\Omega)$ with respect to the norm 
$$|\!|\!|u|\!|\!|_{m}:=\sqrt{  \int_{\R^N} \big(\ln (e+|\xi |)\big)^m |\hat{u}(\xi )|^2   d\xi } $$ 
 and the inner product 
 \begin{align*}
\langle u,v\rangle _m:=    \int_{\R^N} \big(\ln (e+|\xi|)\big)^m \hat{u}(\xi)\hat{v}(\xi)  d\xi \quad{\rm for}\ \, u,v\in \cH_{m,0}(\Omega).
\end{align*}
 
We can have the equivalence of the two spaces. 
 \begin{theorem}\label{th poincare}
For $m\geq 1$,  we have the equivalence of the spaces
$${\bf H}_{m,0}(\Omega)\cong {\cH}_{m,0}(\Omega).$$

 \end{theorem}
 {\bf Proof. }  The critical point is to show the equivalence of the two norms. We claim that for $u\in C_c^1(\Omega)$
  \begin{align}\label{5.1-m}
 |\!|\!|u|\!|\!|_{m}^2\geq \int_{\R^N}(\ln |\xi|)^m |\hat{u}|^2 d\xi \geq  |\!|\!|u|\!|\!|_{m}^2-C_m \|u\|_{L^2}.
\end{align}
It is obvious that  
  \begin{align*}
 \int_{\R^N}\big(\ln (e+|\xi|)\big)^m  |\hat{u}|^2 d\xi \geq \int_{\R^N}(\ln |\xi|)^m |\hat{u}|^2 d\xi. 
\end{align*}
We observe that 
 $$
 \int_{\R^N}(\ln |\xi|)^m |\hat{u}|^2 d\xi  \leq  \int_{\R^N}\big(\ln (e+|\xi|)\big)^m  |\hat{u}|^2 d\xi+  \int_{\R^N}\Big(\big(\ln (e+|\xi|)\big)^m-(\ln |\xi|)^m\Big)  |\hat{u}|^2 d\xi. 
 $$
 Direct computation shows that 
   \begin{align*}
 0 \leq  \int_{\R^N\setminus B_1}\Big(\big(\ln (e+|\xi|)\big)^m-(\ln |\xi|)^m\Big)  |\hat{u}|^2 d\xi  &    \leq 2^m  \int_{\R^N\setminus B_1}\big(\ln (e+|\xi|)\big)^{m-1} \ln (1+\frac{e}{|\xi|}) |\hat{u}|^2 d\xi
 \\&\leq  2^m  \int_{\R^N\setminus B_1}\big(\ln (e+|\xi|)\big)^{m-1} \frac{e}{|\xi|}  |\hat{u}|^2 d\xi
 \\&\leq 2^m \sup_{|\xi|\geq1} \big(\ln (e+|\xi|)\big)^{m-1} \frac{e}{|\xi|}  \, \int_{\R^N\setminus B_1}|\hat{u}|^2 d\xi
 \\&\leq 2^m e \sup_{|\xi|\geq1} \Big(\big(\ln (e+|\xi|)\big)^{m-1} \frac1{|\xi|}\Big)  \, \int_{\R^N}u^2 dx
\end{align*}
 and
    \begin{align*}
0 \leq  \int_{ B_1}\Big(\big(\ln (e+|\xi|)\big)^m-(\ln |\xi|)^m\Big)  |\hat{u}|^2 d\xi
   &    \leq   \ln(1+e) \int_{ B_1}  |\hat{u}|^2 d\xi+ \int_{ B_1} \big(-ln |\xi|)\big)^{m-1}  |\hat{u}|^2 d\xi
 \\&\leq  \ln(1+e)  \|u\|_{L^2(\Omega)}^2  +\||\hat{u}|\|_{L^\infty(\R^N)}^2 \int_{ B_1} \big(-ln |\xi|)\big)^{m-1} d\xi
 \\&\leq \ln(1+e)  \|u\|_{L^2(\Omega)}^2  +\omega_{_N} N^{-m} m!  \|u\|_{L^1(\Omega)}^2  
 \\&\leq \Big(\ln(1+e) +\omega_{_N} N^{-m} m! |\Omega| \Big) \|u\|_{L^2(\Omega)}^2.   
\end{align*}
Thus,  it follows  that for some $c_m>0$ independent of $u$
  \begin{align}\label{5.2-m}
  |\!|\!|u|\!|\!|_{m}^2-  c_m \|u\|_{L^2(\Omega)}^2 \leq \int_{\R^N}(\ln |\xi|)^m |\hat{u}|^2 d\xi  \leq   |\!|\!|u|\!|\!|_{m}^2.  
\end{align}
 
 Furthermore, it follows by (\ref{ineq 1}) with suitable $\epsilon>0$ and $C>0$ that 
    \begin{align}\label{5.5-m}
\frac12 \|u\|^2-C\|u\|_{L^2(\Omega)}^2\leq \cE_m(u,u)  \geq  \frac32 \|u\|^2+C\|u\|_{L^2(\Omega)}^2. 
\end{align}
 Moreover, we see that
\begin{align*}
  |\cA_m(u,u)| &\leq  \sum_{k=1}^m \int_{\R^N} \int_{\R^N}  |u(x)| |u(y)| \, |\j_k(x-y)| dxdy
\\&\leq \sum_{k=1}^m 2^{k-1}\int_{\Omega} \int_{\Omega} |u(x)| |u(y)|  |x-y|^{-N} (\ln |x-y|)^{k-1}1_{B_1^c}(x-y) dxdy
\\&\leq  \sum_{k=1}^m 2^{k-1} r_k^{-N} (\ln r_k)^{k-1}  (\int_{\Omega} |u(x)|dx)^2
\\&\leq  \Big(\sum_{k=1}^m 2^{k-1} r_k^{-N} (\ln r_k)^{k-1}\Big)  |\Omega| \, \|u\|_{L^2(\Omega)}^2
\end{align*}
 where $r_k=e^{\frac{k-1}{N}}$ reaches the maximum of the function  $r\in[1,+\infty)\mapsto r^{-N} (\ln r)^{k-1}$.
 Then we have that 
   \begin{align}\label{5.3-m}
\int_{\R^N}(\ln |\xi|)^m |\hat{u}|^2 d\xi    =  \int_{\R^N} u \cL_m  u d\xi  =  \cE_m(u,u)-\cA_m(u,u) \leq   \frac32 \|u\|^2+C\|u\|_{L^2(\Omega)}^2
\end{align}
and 
   \begin{align}\label{5.4-m}
\int_{\R^N}(\ln |\xi|)^m |\hat{u}|^2 d\xi    =   \cE_m(u,u)-\cA_m(u,u)   \geq  \frac12 \|u\|^2-C\|u\|_{L^2(\Omega)}^2
\end{align}

 As a consequence, we have that 
 \begin{align*}
 |\!|\!|u|\!|\!|_{m}^2-  c_m \|u\|_{L^2(\Omega)}^2   \leq  \int_{\R^N}(\ln |\xi|)^m |\hat{u}|^2 d\xi 
 \leq \frac32 \|u\|_m^2+C\|u\|_{L^2(\Omega)}^2
\end{align*}
 and 
  \begin{align*}
 |\!|\!|u|\!|\!|_{m}^2   \geq  \int_{\R^N}(\ln |\xi|)^m |\hat{u}|^2 d\xi 
 \leq \frac12 \|u\|_m^2-C\|u\|_{L^2(\Omega)}^2, 
\end{align*}
 which implies that for some $c>1$ such that 
  \begin{align*}
\frac1c |\!|\!|u|\!|\!|_{m}    \leq  \|u\|_m  \leq  c  |\!|\!|u|\!|\!|_{m}.
\end{align*} 
 The proof is completed.  \hfill$\Box$\medskip

\noindent{\bf Proof of Theorem \ref{teo 1-m}. }
We first claim that $\cI_m$ is well-defined in $\cH_{m,0}(\Omega) \times \cH_{m,0}(\Omega) $. 

To this end, we have to show $|\cI_m(u,v)|<+\infty$ for any $(u,v)\in \cH_{m,0}(\Omega) \times \cH_{m,0}(\Omega)$. 
 By the definition of $\cI_m$, we only need show that 
$$
 |\cE_m(u,u)|, \ |\cA_m(u,u)|<+\infty    \quad{\rm for}\ \, u\in \cH_{m,0}(\Omega).
$$
by the fact that 
$$|(u,v)_j|\leq  \sqrt{(u,u)_j(v,v)_j}\leq \frac12\big((u,u)_j+(v,v)_j\big).   $$

It follows by (\ref{5.5-m}) that 
$$ |\cE_m(u,u)|= C \|u\|_m^2 <+\infty,$$   
and
$$
  |\cA_m(u,u)|  \leq   C  \int_{\Omega} |u(x)|^2dx <+\infty. 
$$
As a consequence, we have that the functional $\cI_m$ is well-defined in $\cH_{m,0}(\Omega) \times \cH_{m,0}(\Omega) $. \ms

{\it Proof of Part $(i)$.}
The functional $$
  u\in \cH_{m,0}(\Omega)  \mapsto \Phi(u):= \cI_{m}(u,u)
$$
is weakly lower semicontinuous. Recall that 
$$\cP_{m,1}:= \{u\in \cH_{m,0}(\Omega),\,  \norm{u}_{L^2(\Omega)}=1\}.$$
Since the embedding $\cH_{m,0}(\Omega)\hookrightarrow  L^2(\Omega )$ is continuous and compact in Proposition \ref{cr 6.1},
then we have that
$$
 |\Phi(u)|  < +\infty \qquad \text{for $u \in \cP_{m,1}$}
$$
 and $\lambda_{m,1}(\Omega) $ is attained by a function $\phi_1 \in \cP_{m,1}$. Consequently, there exists a Lagrange multiplier $\lambda \in \R$ such that
$$
\cI_{m}(\phi_1, w)=\frac{1}{2}\Phi'(\phi_1)w= \lambda \Big(  \int_{\Omega} \phi_{m,1} w  dx \Big)  \qquad \text{for all $w\in \cH_{m,0}(\Omega)$.}
$$
Choosing $w = \phi_{m,1}$ yields to $\lambda= \lambda_{m,1}(\Omega)$. Hence $\phi_1$ is an eigenfunction of (\ref{eq 1.1}) corresponding to the eigenvalue $\lambda_{m,1}(\Omega)$.   \ms

{\it Proof of Part $(ii)$.} By Theorem \ref{th poincare},  we have that  $\phi_{m,1}\in \cH_{m,0}(\Omega)\cong {\bf H}_{m,0}(\Omega)$, then
$$\lambda_{m,1}(\Omega)  =\inf_{\phi\in  \cH_{m,0}(\Omega)} \cI_{m}(\phi , \phi )=\inf_{\varphi\in {\bf H}_{m,0}(\Omega)}\int_{\R^N}  (\ln |\xi|)^m (\hat{\varphi})^2, $$
then for $m$ even, we have that 
\begin{align*}
\lambda_{m,1}(\Omega)  =  \int_{\R^N} (\ln |\xi|)^m (\hat{\phi}_{m,1})^2 d\xi\geq0. 
\end{align*}
 
{\it Proof of Part $(iii)$.}  Inductively,  we assume that $\phi_{m,1},\dots,\phi_{m,k} \in \cH_{m,0}(\Omega)$ and $\lambda_{m,1}(\Omega)  \le \dots \le \lambda_{m,k}(\Omega)$ are already given for some $k \in \N$ with the properties that for $j=1,\dots,k$, the function $\phi_{m,j}$ is a minimizer
of $\Phi$ within the set
 \begin{align*}
\cP_{m,j}:= \{u\in \cH_{m,0}(\Omega)\::\:  \norm{u}_{L^2(\Omega)}=1,\: \text{$  \int_{\Omega} u\phi_{m,n} \,dx  =0$ for $n=1,\dots j-1$} \},
 \end{align*}
$\displaystyle \lambda_{m,j}(\Omega) = \inf_{\cP_{m,j}} \Phi = \Phi(\phi_{m,j})$, and
\begin{equation}
  \label{eq:inductive-eigenvalue}
  \cI_{m}(\phi_{m,j}, \varphi)=\lambda_{m,j}(\Omega)   \int_{\Omega} \phi_{m,j} \varphi  dx   \qquad \text{for all $\varphi\in \cH_{m,0}(\Omega)$.}
\end{equation}
 
Applying the compact embedding again, we obtain that  the value $\lambda_{m,k+1}(\Omega) $ is attained by a function {$\phi_{m,k+1} \in \cP_{m, k+1}$.} Thus, there exists a Lagrange multiplier $ \lambda= \lambda_{m,k+1}(\Omega)\in \R$ such that
\begin{equation}
  \label{eq:inductive-eigenvalue-k+1}
\cI_{m}(\phi_{m,k+1}, \varphi)=\lambda   \int_{\Omega} \phi_{m,k+1} \varphi  dx   \qquad \text{for all $\varphi \in \cP_{m,k+1}(\Omega)$}
\end{equation}
and by choosing $\varphi = \phi_{m, k+1}$.  Moreover, for $j=1,\dots,k$, we have, by (\ref{eq:inductive-eigenvalue}) and the definition of $\cP_{m, k+1}(\Omega)$, that
\begin{align*}
\cI_{m}(\phi_{m,k+1}, \phi_{m,j})&=\cI_{m}(\phi_{m,j}, \phi_{m,k+1}) 
\\&= \lambda_{m,j}(\Omega)  \int_{\Omega} \phi_{m,j} \phi_{m, k+1} \,dx  
 \\& = 0
\\&= \lambda_{m,k+1}(\Omega)  \int_{\Omega} \phi_{m,j} \phi_{m, k+1} \,dx.
\end{align*}
Hence (\ref{eq:inductive-eigenvalue-k+1}) holds with $\lambda=  \lambda_{m,k+1}(\Omega) $ for all $\varphi\in \cH_{m,0}(\Omega)$.
Inductively, we have now constructed a normalized sequence $(\phi_{m,k})_{k\in\N}$ in $\cH_{m,0}(\Omega)$ and a nondecreasing sequence $\{\lambda_{m,k}(\Omega) \}_{k\in\N}$ in $\R$ such that $\phi_k$ is an eigenfunction of (\ref{eq 1.1}) corresponding to $\lambda =  \lambda_{m,k}(\Omega)  $ for every $k \in \N$. Moreover, by construction, the sequence $\{\phi_k\}_k$ forms an orthonormal system in $L^2(\Omega)$.  

{\it We now show  $\lim \limits_{k\to+\infty} \lambda_{m,k}(\Omega)  =+\infty.$}  Assume by contradiction that $c:= \lim \limits_{k \to \infty}\lambda_{m,k}(\Omega) <+\infty$, then we deduce that
$ (\phi_{m,k},\phi_{m,k})_m\le c$ for every $k \in \N$. Hence the sequence $(\phi_{m,k})_{k\in \N}$ is bounded in $\cH_{m, 0}(\Omega)$, and therefore by Rellich compactness theorem, $(\phi_{m,k})_k$
contains a convergent subsequence $(\phi_{m,k_j})_j$ in $L^2(\Omega)$. However, this   is impossible since the functions $\{\phi_{m,k_j}\}_{j \in \N}$ are $L^2(\Omega)$-orthogonal.\ms 

 {\it Proof of Part $(iv)$.}  We prove that $\{\phi_{m,k}\::\: k \in \N\}$ is an orthonormal basis of $L^2(\Omega)$.    By contradiction, we assume that there exists $v \in \cH_{m, 0}(\Omega)$ with $\|v\|_{L^2(\Omega)} = 1$
 and $ \int_{\Omega} v \phi_k\,dx   =0$ for any $k \in \N$.
Since $\lim \limits_{k\to\infty} \lambda_{m,k}(\Omega)=+\infty$, there exists an integer $k_0>0$  such that
$$\Phi(v)<\tilde \lambda_{m,k_0} =\inf_{u\in \cP_{m,k_0}}\Phi(u),$$
 which, by definition of $\cP_{m, k_0}$,  implies that
$\int_{\Omega} v \phi_{m,k}\,dx  \not = 0$ for some $k \in \{1,\dots,k_0-1\}$, which leads to a contradiction. 

We conclude that $\cH_{m, 0}(\Omega)$ is contained in the $L^2(\Omega)$-closure of the span of $\big\{\phi_{m,k}:\: k \in \N\big\}$. Since $\cH_{m, 0}(\Omega)$ is dense in $L^2(\Omega)$, we conclude that the span of $\big\{\phi_{m,k}:\: k \in \N\big\}$ is dense in $L^2(\Omega)$, and hence $\big\{\phi_{m,k}:\: k \in \N\big\}$ is an orthonormal basis of $L^2(\Omega)$.\ms 

 {\it Proof of Part $(v)$.}
  Let 
\begin{align}\label{hh}
{\bf h}_m(t):= \sum^m_{j=1} \alpha_j  (-2\ln t)^{j-1} 1_{(0,r_0)}(t)\quad {\rm for}\ t>0.
\end{align}
 Since $\alpha_m=m \kappa_1(0)>0$, for suitable $r_0\in(0,\frac12]$ such that   
 $$\sum^m_{j=0} \alpha_j  (-2\ln t)^{j-1} >0\quad {\rm for}\ 0<t<r_0 $$

 For $\Omega \subset B_{r_0}$  and $\alpha_m>0$, then for all $u\in \cH_{m,1}(\Omega)$,   $\cE_m$ reduces to
\begin{align*}
 \cQ_m (u, u)=\alpha_0\int_\Omega u^2 dx+ \int_{\R^N}\int_{\R^N} \frac{\big(u(x)-u(y)\big)^2    }{|x-y|^{N}} {\bf h}_m(|x-y|)    dxdy>0 
\end{align*}
and 
\begin{align*}
  \cQ_m (|u|, |u|) &=  \alpha_0\int_\Omega |u|^2 dx+ \int_{\R^N}  \int_{\R^N} \frac{\big(|u(x)|-|u(y)|\big)^2    }{|x-y|^{N}}  {\bf h}_m(|x-y|)dxdy  
\\& \ \leq \alpha_0\int_\Omega u^2 dx+ \int_{\R^N}  \int_{\R^N} \frac{\big(u(x)-u(y)\big)^2    }{|x-y|^{N}}   {\bf h}_m(|x-y|)  dxdy\ = \cQ_m (u, u),
\end{align*}
where the equality holds when $u$ does not change sign.
 
   Let $w_0 \in \cH_{m,1}(\Omega)$ be a $L^2$-normalized eigenfunction of $\cL_m$ corresponding to the eigenvalue $\lambda_{m,1}(\Omega)$, i.e. we have
\begin{equation}
  \label{eigenvalue-k=1}
 \cQ_m(w_0, \phi)= \lambda_{m,1}(\Omega)  \int_\Omega w_0 \phi\, dx \qquad \text{for all $\phi \in \mathbb{H}_{m,1}(\Omega)$.}
\end{equation}

We show that $w_0$ does not change sign. Indeed, choosing $\phi= w_0$ in (\ref{eigenvalue-k=1}), we see that $w$ is a minimizer of $\cE_m(w,w)\big|_{w\in \cP_{m,1}}$. On the other hand, we also have $|w| \in \cP_{m,1}$ and
$$
\lambda_{m,1}(\Omega)=\inf_{\phi \in \cH_{m,0}(\Omega)} \cQ_m(\phi,\phi)  \leq    \cQ_m(|w_0|,|w_0|) \leq \cQ_m(w_0,w_0) =\lambda_{m,1}(\Omega),
$$
 then  $w_0$ must not be sign-changing. In particular, we may assume that $w_0$ is nonnegative.  Observe that
$$
({\bf Q}_m w_0, \phi )_0 -\lambda_{m,1}(\Omega)(w_0, \phi )_0=  0 \qquad \text{for all $\phi \in C^\infty_c(\Omega)$, $\phi \ge 0$,}
$$
where $${\bf Q}_m u(x)= \int_{\R^N}\frac{ u(x)-u(y)}{|x-y|^N}  {\bf h}_m(x-y)dy+\alpha_0 u(x).$$
  Hence $w_0$ is a nontrivial, nonnegative weak  solution of the equation 
$${\bf Q}_m w_0 - \lambda_{m,1}(\Omega)  w_0=0\quad{\rm in}\  \Omega.$$
 Therefore, we apply \cite[Theorem 1.1]{JW-preprint}  to yield that $w_0>0$ a.e. in $\Omega$.\vs

Now suppose by contradiction that there is a function $w \in \cP_{m,1}(\Omega)$ satisfying (\ref{eigenvalue-k=1}) and such that $w \not = t\phi_{m,1}$ for every $t \not=\pm1$. Then there exist a linear combination $\tilde w$ of $w$ and $\phi_{m,1}$ which changes sign, and we may also assume that $\tilde w$ is $L^2$-normalized. Since $\tilde w$ also satisfies (\ref{eigenvalue-k=1}) in place of $w$, we arrive at a contradiction.  Since the eigenvalue $0<\lambda_{m,1}(\Omega)$ is simple, then  $0<\lambda_{m,1}(\Omega)< \lambda_{m,2}(\Omega)$.
\hfill$\Box$

\subsection{Faber-Krahn inequality}

To prove Theorem \ref{sec:faber-Krahn}, we need following lemma, which is from \cite[Theorem 3]{B}
in $\R^N$.

\begin{lemma}  \label{lm 2.1} Assume that  
\begin{enumerate}
\item[(a)] $\Psi:[0,+\infty)\to[0,+\infty)$ is a positive convex, monotone increasing  function such that $\Psi(0)=0$ 
  $\Psi''\geq 0$ and  $t\Psi'(t)$ convex; 
  \item[(b)] $K:(0,+\infty)\to[0,+\infty)$ is monotone decreasing; 
   \item[(c)] $\rho:[0,+\infty)\to[0,+\infty)$ is monotone increasing.
\end{enumerate}
 Let $u^*$ be the equi-measurable rearrangement of $u$.
 
Then for $f,g$ nonnegative functions, 
 $$\int_{\R^N}\int_{\R^N} \Psi\Big[\frac{f(x)-g(y))}{\rho(|x-y|)} \Big]K(|x-y|)dxdy\geq    \int_{\R^N}\int_{\R^N} \Psi\Big[\frac{f*(x)-g^*(y))}{\rho(|x-y|)}\Big] K(|x-y|)dxdy. $$
 
Moreover, if $K$ is strictly decreasing and $\Psi$ is strictly convex, 
the above equality holds only when $f=f^*$ and $g=g^*$
 
\end{lemma}

\medskip

\noindent{\bf Proof of Theorem \ref{sec:faber-Krahn}. }  We assume that $\Omega\subset B_{\bar r}$ with $\bar r<\frac12 r_m$.  Recall that ${\bf h}_{m}$ is strictly decreasing in $(0,r_m)$.
  The Rayleigh-Ritz principle characterizes the principal eigenvalue as a minimum
\[
\lambda_{m,1}(\Omega)=\inf_{u\in \cH_{m,0}(\Omega)\setminus\{0\}}\frac{\cQ_m(u)}{\|u\|^2_{L^2(\Omega)}}, 
 \]
 where 
 $$
\cQ_m(u)= \alpha_0 \int_{\Omega} u^2 dx+   \int_{\R^N}\int_{\R^N}\frac{\big(u(x)-u(y)\big)^2}{|x-y|^{N}}{\bf h}_m(|x-y|)\, dxdy. $$
 Given $\epsilon\in(0,\frac14 \bar r)$, we let ${\bf h}_{m,0}: [0,+\infty)\to (0,+\infty)$ be  a strictly  decreasing continuous function  such that 
 $${\bf h}_{m,0}(t)  =\left \{
  \begin{array}{lll}
  {\bf h}_{m}(t)  &\quad  {\rm for} \ t\in (0, 2\bar r),\\[2mm]
   {\bf h}_{m}(2\bar r) e^{2\bar r-t} &\quad {\rm if} \ s\in [2\bar r,+\infty).
  \end{array}
\right.
  $$

Let $u\in \cH_{m,0}(\Omega)$ be a nonnegative function,   we apply Lemma \ref{lm 2.1} with $\Psi(t)=t^2$, $\rho(t)=t^{\frac N2}$, 
 $K(t)= {\bf h}_{m,0}(t)$ for  $t\in[0,+\infty)$, which is strictly decreasing,   and $f=g=u$,  to obtain that 
  \begin{align*}
 \int_{\R^N}\int_{\R^N}\frac{\big(u^*(x)-u^*(y)\big)^2}{|x-y|^{N }}{\bf h}_{m,0}(|x-y|) \, dxdy  \leq       \int_{\R^N}\int_{\R^N}\frac{\big(u(x)-u(y)\big)^2}{|x-y|^{N}}{\bf h}_{m,0}(|x-y|) \, dxdy,
    \end{align*}
   where the equality holds only when $u=u^*$.

If $\Omega\not=B_{\bar r}$,  then 
 \begin{align*}
 \int_{\R^N}\int_{\R^N}   \frac{\big(u^*(x)-u^*(y)\big)^2}{|x-y|^{N }} {\bf h}_{m,0}(|x-y|) \, dxdy  <      \int_{\R^N}\int_{\R^N}\frac{\big(u(x)-u(y)\big)^2}{|x-y|^{N}}{\bf h}_{m,0}(|x-y|) \, dxdy,
    \end{align*}
then 
 \begin{align*}
  \int_{\R^N}\int_{\R^N}\frac{\big(u^*(x)-u^*(y)\big)^2}{|x-y|^{N }}{\bf h}_m(|x-y|) \, dxdy  
  & = \int_{ B_{r_m}} \int_{B_{r_m}}\frac{\big(u^*(x)-u^*(y)\big)^2}{|x-y|^{N}}{\bf h}_{m,0}(|x-y|)\, dxdy  
 \\&=  \int_{\R^N}\int_{\R^N}\frac{\big(u^*(x)-u^*(y)\big)^2}{|x-y|^{N}}{\bf h}_{m,0}(|x-y|)\, dxdy
  \\& <\int_{\R^N}\int_{\R^N}\frac{\big(u(x)-u(y)\big)^2}{|x-y|^{N}}{\bf h}_{m,0}(|x-y|)\, dxdy
\\&=  \int_{\R^N}\int_{\R^N}\frac{\big(u(x)-u(y)\big)^2}{|x-y|^{N}}{\bf h}_m(|x-y|) \, dxdy. 
    \end{align*}

We recall  that   for any $f\in L^{p}(\Omega)$ nonnegative with $p>1$, 
 $$\int_{B_r}(f^{*})^p dx= \int_{\Omega}f^pdx. $$

 As a consequence,  if $\Omega\not=B_{\bar r}$ we have that 
  \begin{align*}
\lambda_{m,1}(\Omega) &=  \frac{\alpha_0\int_{\Omega} \phi_{m,1}^2 dx+   \int_{\R^N}\int_{\R^N}\frac{(\phi_{m,1}(x)-\phi_{m,1}(y))^2}{|x-y|^{N}}{\bf h}_m(|x-y|)\, dxdy }{ \int_{\Omega}   \phi_{m,1} ^2 dx }   
 \\[2mm]&> \frac{\alpha_0 \int_{ B_r} (\phi_{m,1}^*)^2 dx+   \int_{\R^N}\int_{\R^N}\frac{(\phi_{m,1}^*(x)-\phi_{m,1}^*(y))^2}{|x-y|^{N }}{\bf h}_m(|x-y|)\, dxdy }{ \int_{B_{\bar r}}   (\phi_{m,1}^*) ^2 dx} 
  \\[2mm]& \geq \inf_{\substack{u\in \cH_{m,0}(B_{\bar r}) \setminus\{0\} }}\frac{\cQ_m(u,u)}{\|u\|^2_{L^2(B_{\bar r})}}
   \\[2mm]& =\lambda_{m,1}(B_{\bar r})
\end{align*}
as claimed. 

Moreover, when $\Omega=B_{\bar r}$, we see that  
  \begin{align*}
\lambda_{m,1}(B_{\bar r})  &= \frac{\alpha_0\int_{B_{\bar r}} \phi_{m,1}^2 dx+   \int_{\R^N}\int_{\R^N}\frac{(\phi_{m,1}(x)-\phi_{m,1}(y))^2}{|x-y|^{N}}{\bf h}_m(|x-y|)\, dxdy }{ \int_{B_{\bar r}}   \phi_{s,1} ^2 dx }   
 \\[2mm]&\geq \frac{\alpha_0 \int_{ B_{\bar r}} (\phi_{m,1}^*)^2 dx+   \int_{\R^N}\int_{\R^N}\frac{(\phi_{m,1}^*(x)-\phi_{m,1}^*(y))^2}{|x-y|^{N+2s}}{\bf h}_m(|x-y|)\, dxdy }{ \int_{B_{\bar r}}   (\phi_{m,1}^*) ^2 dx } 
  \\[2mm]& \geq \inf_{\substack{u\in \cH_{m,s}(B_{\bar r})\setminus\{0\}}}\frac{\cQ_m(u,u)}{\|u\|^2_{L^2(B_{\bar r})}},
\end{align*}
which, together with the uniqueness, implies that $\phi_{m,1}=\phi_{m,1}^*$, that is, $\phi_{m,1}$ is radially symmetric and decreasing with respect to $|x|$. 
  \hfill$\Box$\medskip

  \setcounter{equation}{0}
 \appendix
\section{Appendix }
\subsection{ Dirichlet eigenvalues for $\cK_m$}

We consider the bounds of the Dirichlet eigenvalues of $\cK_m$
  \begin{equation}\label{eq 1.1-mp}
\left\{ \arraycolsep=1pt
\begin{array}{lll}
\cK_m  u=\mu u\quad \  &{\rm in}\ \,   \Omega,\\[1.5mm]
 \phantom{   \cL_m    }
  u=0    &{\rm in}\ \, \R^N\setminus \Omega, 
\end{array}
\right.
\end{equation}
where $\Omega$ is a bounded Lipschitz domain in $\R^N$. 

A function $u\in \cH_{m,0}(\Omega)$ is called the eigenfunction of \eqref{eq 1.1-mp} corresponding to the eigenvalue  $\lambda$ if
\begin{equation}\label{weak-m}
\cG_m(u,\phi)= \mu\int_{\Omega}u\phi dx,\quad \forall  \, \phi\in \cH_{m,0}(\Omega),
\end{equation}
where  
  \begin{align*}
\cG_m(u,\phi) &= (u,\phi)_m-(\,\j_m\ast u, \phi)_0
 \\[1mm]& = \int_{\R^N} \int_{\R^N}\big(u(x)-u(y)\big)\big(\phi(x)-\phi(y)\big) \k_m(x-y)dxdy-\int_{\R^N} \int_{\R^N} u(x)\phi(y) \,\j_m(x-y)dxdy.
\end{align*}
Note that $\cG_m$ is a closed, symmetric and semi-bounded quadratic form with domain $\cH_{m,0}(\Omega)\subset L^2(\Omega)$.

     \begin{theorem}\label{teo 1-mm}
Assume that  $\Omega\subset \R^N$ be a bounded  Lipschitz  domain.  

\begin{itemize}
\item[(i)] Problem \eqref{eq 1.1-mp}  admits  an eigenvalue $\mu_{m,1}(\Omega)>0$   characterized  by
\begin{equation}\label{Lambda-1-s-m}
\mu_{m,1}(\Omega)=\inf_{\substack{u\in \cG_{m,0}(\Omega)\\ u\neq 0}}\frac{\cG_m(u,u)}{\|u\|^2_{L^2(\Omega)}}= \inf_{u\in \cM_{m,1}(\Omega)}\cG_m(u,u),
\end{equation}
with $$\cM_{m,1}(\Omega):=\{u\in \cH_{m,0}(\Omega): \|u\|_{L^2(\Omega)}=1\}$$
 and there exists a positive function $\varphi_{m,1}\in \cM_{m,1}(\Omega)$ achieving the minimum in \eqref{Lambda-1-s-m}, 
i.e. 
 $
 \mu_{m,1}(\Omega)= \cG_{m}(\varphi_{m,1},\varphi_{m,1}).$ 
 
\item[(ii)]  The first eigenvalue $\mu_{m,1}(\Omega)$ is simple, that is, if $u\in  \cH_{m,s}(\Omega)$ satisfies \eqref{weak-m}  with $\lambda=\mu_{m,1}(\Omega)$,  then $u=\alpha\varphi_{m,1}$ for some $\alpha\in \R$.

\item[(iii)] Problem \eqref{eq 1.1-mp}  admits a sequence of eigenvalues $\{\mu_{m,k}(\Omega)\}_{k\in \N}$ such that 
\[
 \mu_{m,1}(\Omega)< \mu_{m,1}(\Omega)\le \cdots\le \mu_{m,k}(\Omega)\le \mu_{m,k+1}(\Omega)\leq \cdots,
\]
with corresponding eigenfunctions $\varphi_{m,k}$ with $k\in \N$ and~
$
\displaystyle \lim_{k\to \infty}\mu_{m,k}(\Omega) = +\infty.
$

Moreover, for any $k\in\N$, the  eigenvalue $\mu_{m,k}(\Omega)$ can be characterized as
\begin{equation}\label{Lambda-k-s}
\mu_{m,k}(\Omega)= \inf_{u\in \cM_{m,k}(\Omega)}\cG_m(u,u),
\end{equation}
where 
\[
\cM_{m,k}(\Omega):= \Big\{ u\in \cH_{m,0}(\Omega): \int_{\Omega}u\varphi_{m,j}dx  =0 \ \text{ for \ } j=1,2,\cdots k-1, \  \ \|u\|_{L^2(\Omega)}=1\Big\}.
\]
\item[(iv)] The sequence of eigenfunctions  $\{\varphi_{m,k}\}_{k\in\N}$ corresponding to eigenvalues $\{\mu_{m,k}(\Omega)\}$ form a complete orthonormal basis of $L^2(\Omega)$ and  an orthogonal system of $\cH_{m,0}(\Omega)$.

\item[(v)]  we assume more that  $\Omega\subset B_{\frac12}$, then $\mu_{m,1}(\Omega)>0$. 
\end{itemize}
 \end{theorem} 
 
     \begin{lemma}
\label{fram-1-modulus}
Let $m\geq 1$ $u \in \cH_{m,0}(\Omega)$, then  $|u| \in \cH_{m,0}(\Omega)$  and
\begin{equation}
  \label{eq:modulus-invariance}
\cG_m(|u|,|u|)\le \cG_m(u,u).
\end{equation}
Moreover, equality holds in (\ref{eq:modulus-invariance}) if and only if $u$ does not change sign.
\end{lemma}
 {\bf Proof. }    We observe that 
\begin{align}\label{abs-a1}
 (|u|,|u|)_m \leq   (u,u)_m,
\end{align}
where the equality holds only when $u$ does not change sign. 
 Moreover, we have that 
 \begin{align*}
   \int_{\R^N}  \int_{\R^N} |u(x)| |u(y)|  \,\j_m(x-y) dxdy \geq    \int_{\R^N}  \int_{\R^N}  u(x) u(y)   \,\j_m(x-y) dxdy,
\end{align*}
which, together with (\ref{abs-a1}),  implies that  $\cG_m(|u|,|u|)\le \cG_m(u,u)$ and its equality holds when $u$ does not change sign.\hfill$\Box$\medskip

\noindent{\bf Proof of Theorem \ref{teo 1-mm}. }     The proof are similar the one of Theorem \ref{teo 1-m} replaced $\cE_m$ by $\cG_m$.

 Here we would like mention that  in general bounded domain and $w$ is an eigenfunction corresponding to $\mu_{m,1}(\Omega)>0$,   we also have $|w| \in \cM_{m,1}$ and
$$
\mu_{m,1}(\Omega)=\inf_{\phi \in \cH_{m,0}(\Omega)} \cG_m(\phi,\phi)  \leq    \cG_m(|w|,|w|) \leq \cG_m(w,w) = \mu_{m,1}(\Omega),
$$
by Lemma~\ref{fram-1-modulus} and   $w$ does not change sign.
     
    Moreover,  we note that   $\varphi_{m,1}=0$ in $\R^N\setminus \Omega$ and $\Omega\subset B_{\frac12}$, then
   \begin{align*}
   \int_{\Omega} [\,\j_m * \varphi_{m,1}] \varphi_{m,1}\,dx=\int_{|x-y|>1} \frac{u(x)u(y)}{|x-y|^N} (-2\ln |x-y|)^{m-1}dxdy =0,  
  \end{align*}
 which implies that $\mu_{m,1}(\Omega)>0$.     \hfill$\Box$  \ms
 
 From the proof of Theorem \ref{sec:faber-Krahn},  we have  the  Faber-Krahn inequality of the first eigenvalue for the $\cK_m$.  
 \begin{theorem}  
 Let $\Omega\subset B_{\frac12}$ be a Lipschitz domain and   $|\Omega| =|B_{\bar r}|$ for some $\bar r\in(0,\frac12]$,  then only the ball $B_{\bar r}$  minimizes $\mu_{m,1}$, i.e.
$$\mu_{m,1}(\Omega)\geq \mu_{m,1}(B_{\bar r}),$$
where $'='$ holds only for $\Omega=B_{\bar r}$. 
Moreover, the first eigenfunction $\varphi_{m,1}$ for $\Omega=B_{\bar r}$ is radially symmetric and decreasing with respective to $|x|$.  

\end{theorem}


  \bigskip
  
 \medskip 
  
 { \small 
  
  \noindent {\bf Acknowledgements: }   The author is grateful to   Professor Tobias Weth for valuable discussions.    \medskip
  
  \noindent {\bf Funds: }   H. Chen is supported by  NNSF of China (Nos. 12071189, 12361043) 
and Jiangxi Natural Science Foundation (No. 20232ACB201001). 
  
  }
  
  \medskip
 

\end{document}